\documentclass{amsart}
\def\MMM{}
\def\WW{}
\def\FORM{}

\numberwithin{equation}{section}
\numberwithin{figure}{section}
\usepackage{amssymb}
\let\cal\mathcal

\def\Bscr{{\cal B}}
\def\Cscr{{\cal C}}
\def\Dscr{{\cal D}}
\def\Escr{{\cal E}}

\def\Lscr{{\cal L}}

\def\Oscr{{\cal O}}
\def\Pscr{{\cal P}}

\def\Sscr{{\cal S}}

\def\Uscr{{\cal U}}
\def\Vscr{{\cal V}}
\def\Wscr{{\cal W}}

\newcommand{\proj}{\operatorname{proj}}

\def\pr{\mathop{\text{pr}}\nolimits}

\def\Ab{\mathbb{Ab}}

\def\Mod{\operatorname{Mod}}
\def\mod{\operatorname{mod}}

\def\Sk{\operatorname{Sk}}
\def\gr{\operatorname{gr}}

\def\Ch{\mathop{\mathrm{Ch}}}

\def\gr{\operatorname {gr}}
\def\Spec{\operatorname {Spec}}

\def\Ext{\operatorname {Ext}}
\def\Hom{\operatorname {Hom}}

\def\ker{\operatorname {ker}}

\def\Tor{\operatorname {Tor}}

\def\r{\longrightarrow}

\DeclareMathOperator{\Ind}{Ind}

\DeclareMathOperator{\Aut}{Aut}

\newtheorem{lemma}{Lemma}[section]
\newtheorem{proposition}[lemma]{Proposition}
\newtheorem{theorem}[lemma]{Theorem}
\newtheorem{corollary}[lemma]{Corollary}

\theoremstyle{definition}

\newtheorem{example}[lemma]{Example}
\newtheorem{definition}[lemma]{Definition}

{

}

\theoremstyle{remark}

\newtheorem{remark}[lemma]{Remark}

\newdimen\uboxsep \uboxsep=1ex
\def\uboxn#1{\vtop to 0pt{\hrule height 0pt depth 0pt\vskip\uboxsep
\hbox to 0pt{\hss #1\hss}\vss}}

\def\uboxs#1{\vbox to 0pt{\vss\hbox to 0pt{\hss #1\hss}
\vskip\uboxsep\hrule height 0pt depth 0pt}}

\usepackage{amscd}
\usepackage[all]{xy}

\def\redef#1#2{\def#1{#2}}
\redef{\ovl}{\overline}
\redef{\Kern}{\mathrm{ker}}
\redef{\Ob}{\mathrm{Ob}}
\redef{\Mor}{\mathrm{Mor}}
\redef{\Ch}{\mathrm{Ch}}
\redef{\K}{\mathrm{K}}
\redef{\D}{\mathrm{D}}
\redef{\Cokern}{\mathrm{Coker}}
\redef{\Beeld}{\mathrm{Im}}
\redef{\Sub}{\mathsf{Sub}}
\renewcommand{\lim}{\mathrm{lim}}
\redef{\colim}{\mathrm{colim}}
\redef{\lrank}{\mathrm{rank}}
\redef{\cham}{\mathrm{cham}}
\redef{\codim}{\mathrm{codim}}
\redef{\Aut}{\mathrm{Aut}}
\redef{\proj}{\mathrm{proj}}
\redef{\flag}{\mathrm{flag}}
\redef{\type}{\mathrm{type}}
\redef{\Tor}{\mathrm{Tor}}
\redef{\Ext}{\mathrm{Ext}}
\redef{\Hom}{\mathrm{Hom}}
\redef{\Def}{\mathrm{Def}}
\redef{\ddef}{\mathrm{def}}
\redef{\deff}{\mathrm{def}}
\redef{\Sk}{\mathrm{Sk}}
\redef{\pr}{\mathrm{pr}}
\redef{\inn}{\mathrm{in}}
\redef{\es}{\mathrm{es}}
\redef{\op}{^{\mathrm{op}}}
\redef{\Z}{\mathbb{Z}}
\redef{\Q}{\mathbb{Q}}
\redef{\R}{\mathbb{R}}
\redef{\N}{\mathbb{N}}
\redef{\C}{\mathbb{C}}
\redef{\AAA}{\mathfrak{a}}
\redef{\BBB}{\mathfrak{b}}
\redef{\CCC}{\mathfrak{c}}
\redef{\DDD}{\mathfrak{d}}
\redef{\III}{\mathfrak{i}}
\redef{\JJJ}{\mathfrak{j}}
\redef{\PPP}{\mathfrak{p}}
\redef{\SSS}{\mathfrak{s}}
\redef{\TTT}{\mathfrak{t}}
\redef{\XXX}{\mathfrak{x}}
\redef{\UUU}{\mathfrak{u}}
\redef{\VVV}{\mathfrak{v}}
\redef{\WWW}{\mathfrak{w}}
\redef{\GGGG}{\mathfrak{g}}
\redef{\LLLL}{\mathfrak{l}}
\redef{\Set}{\ensuremath{\mathsf{Set}} }
\redef{\Ab}{\ensuremath{\mathsf{Ab}} }
\redef{\Mod}{\ensuremath{\mathsf{Mod}} }
\redef{\PMod}{\ensuremath{\mathsf{PreMod}} }
\redef{\Pre}{\ensuremath{\mathsf{Pr}} }
\redef{\Sh}{\ensuremath{\mathsf{Sh}} }
\redef{\mmod}{\ensuremath{\mathsf{mod}} }
\redef{\RMod}{\ensuremath{\mathsf{Mod}(R)} }
\redef{\SMod}{\ensuremath{\mathsf{Mod}(S)} }
\redef{\Rmod}{\ensuremath{\mathsf{mod}(R)} }
\redef{\Smod}{\ensuremath{\mathsf{mod}(S)} }
\redef{\ModR}{\ensuremath{\mathsf{Mod}(R)} }
\redef{\ModS}{\ensuremath{\mathsf{Mod}(S)} }
\redef{\modR}{\ensuremath{\mathsf{mod}(R)} }
\redef{\modS}{\ensuremath{\mathsf{mod}(S)} }

\redef{\AAAMod}{\ensuremath{{\mathsf{Mod}}{(\mathfrak{a})}} }
\redef{\ModAAA}{\ensuremath{{\mathsf{Mod}}{(\mathfrak{a})}} }
\redef{\BBBMod}{\ensuremath{{\mathsf{Mod}}{(\mathfrak{b})}} }
\redef{\ModBBB}{\ensuremath{{\mathsf{Mod}}{(\mathfrak{b})}} }
\redef{\AAAmod}{\ensuremath{{\mathsf{mod}}{(\mathfrak{a})}} }
\redef{\modAAA}{\ensuremath{{\mathsf{mod}}{(\mathfrak{a})}} }
\redef{\BBBmod}{\ensuremath{{\mathsf{mod}}{(\mathfrak{b})}} }
\redef{\modBBB}{\ensuremath{{\mathsf{mod}}{(\mathfrak{b})}} }

\redef{\RCat}{\ensuremath{{\mathsf{Cat}}(R)} }
\redef{\SCat}{\ensuremath{{\mathsf{Cat}}(S)} }
\redef{\CatR}{\ensuremath{{\mathsf{Cat}}(R)} }
\redef{\CatS}{\ensuremath{{\mathsf{Cat}}(S)} }

\redef{\Rng}{\ensuremath{\mathsf{Rng}} }
\redef{\Cat}{\ensuremath{\mathsf{Cat}} }
\redef{\PreCat}{\ensuremath{\mathsf{PreCat}} }
\redef{\Gd}{\ensuremath{\mathsf{Gd}} }
\redef{\PreCX}{\ensuremath{\mathsf{Pre}_{\mathcal{C}}(X)}}
\redef{\PreCY}{\ensuremath{\mathsf{Pre}_{\mathcal{C}}(Y)}}
\redef{\PreoX}{\ensuremath{\mathsf{Pre}(\theta _X)}}
\redef{\PreAbX}{\ensuremath{\mathsf{Pre}_{\Ab}(X)}}
\redef{\PreRngX}{\ensuremath{\mathsf{Pre}_{\Rng}(X)}}
\redef{\SchCX}{\ensuremath{\mathsf{Sch}_{\mathcal{C}}(X)}}
\redef{\SchCY}{\ensuremath{\mathsf{Sch}_{\mathcal{C}}(Y)}}
\redef{\SchAbX}{\ensuremath{\mathsf{Sch}_{\Ab}(X)}}
\redef{\SchRngX}{\ensuremath{\mathsf{Sch}_{\Rng}(X)}}
\redef{\CatX}{\ensuremath{\mathsf{Cat} (X)}}
\redef{\CatY}{\ensuremath{\mathsf{Cat} (Y)}}
\redef{\CatXp}{\ensuremath{\mathsf{Cat} (X,p)}}
\redef{\CatXA}{\ensuremath{\mathsf{Cat} (X,A)}}
\redef{\CatXfU}{\ensuremath{\mathsf{Cat} (X,f(U))}}
\redef{\PreC}{\ensuremath{\mathsf{Pr}(\mathcal{C})}}
\redef{\IndC}{\ensuremath{\mathsf{Ind}(\mathcal{C})}}
\redef{\Ind}{\ensuremath{\mathsf{Ind}}}
\redef{\Fun}{\ensuremath{\mathsf{Fun}}}
\redef{\Add}{\ensuremath{\mathsf{Add}}}
\redef{\Inj}{\ensuremath{\mathsf{Inj}}}
\redef{\Cof}{\ensuremath{\mathsf{Cof}}}
\redef{\inj}{\ensuremath{\mathrm{Inj}}}
\redef{\cof}{\ensuremath{\mathrm{Cof}}}
\redef{\Fp}{\ensuremath{\mathsf{Fp}}}
\redef{\ShCT}{\ensuremath{\mathsf{Sh}(\mathcal{C},\mathcal{T})}}
\redef{\ShCL}{\ensuremath{\mathsf{Sh}(\mathcal{C},\mathcal{L})}}
\redef{\ShCLnul}{\ensuremath{\mathsf{Sh}(\mathcal{C},\mathcal{L}
_0)}}
\redef{\ShCLplus}{\ensuremath{\mathsf{Sh}(\mathcal{C},\mathcal{L}
^{+})}}
\redef{\SepCT}{\ensuremath{\mathsf{Sep}(\mathcal{C},\mathcal{T})}}
\redef{\SepCL}{\ensuremath{\mathsf{Sep}(\mathcal{C},\mathcal{L})}}
\redef{\SepCLnul}{\ensuremath{\mathsf{Sep}(\mathcal{C},\mathcal{L}
_0)}}
\redef{\SepCLplus}{\ensuremath{\mathsf{Sep}(\mathcal{C},\mathcal{L}
^{+})}}
\redef{\lra}{\longrightarrow}
\redef{\ra}{\rightarrow}
\redef{\Stensor}{\ensuremath{(S \otimes_R -)}}
\redef{\HomS}{\ensuremath{\mathrm{Hom}_R(S,-)}}

\newcommand{\aaa}{\ensuremath{\mathcal{A}}}
\newcommand{\bbb}{\ensuremath{\mathcal{B}}}
\newcommand{\ccc}{\ensuremath{\mathcal{C}}}
\newcommand{\ddd}{\ensuremath{\mathcal{D}}}
\newcommand{\eee}{\ensuremath{\mathcal{E}}}

\newcommand{\iii}{\ensuremath{\mathcal{I}}}

\newcommand{\kkk}{\ensuremath{\mathcal{K}}}
\newcommand{\LLL}{\ensuremath{\mathcal{L}}}

\newcommand{\ooo}{\ensuremath{\mathcal{O}}}

\newcommand{\rrr}{\ensuremath{\mathcal{R}}}
\newcommand{\sss}{\ensuremath{\mathcal{S}}}

\newcommand{\uuu}{\ensuremath{\mathcal{U}}}
\newcommand{\vvv}{\ensuremath{\mathcal{V}}}
\newcommand{\www}{\ensuremath{\mathcal{W}}}

\title{Deformation theory of abelian categories}
\author{Wendy T. Lowen}
\address{Departement DWIS\\ Vrije Universiteit Brussel\\ Pleinlaan
2\\1050 Brussel\\ Belgium}
\email[Wendy T. Lowen]{wlowen@vub.ac.be}
\author{Michel Van den Bergh}
\address{Departement WNI\\Limburgs Universitair
Centrum\\ Universitaire Campus\\ Building D\\ 3590
Diepenbeek\\ Belgium}
\email[Michel Van den Bergh]{vdbergh@luc.ac.be}
\thanks{The first author is an aspirant at the FWO}
\thanks{The second
author is a senior researcher at the FWO}
\begin{document}
\begin{abstract}
     In this paper we develop the basic infinitesimal deformation theory
     of \emph{abelian categories}.  This theory yields a natural
     generalization of the well-known deformation theory of algebras
     developed by Gerstenhaber.  As part of our deformation theory we
     define a notion of flatness for abelian categories. We show that
     various basic properties are preserved under flat deformations and
     we construct several equivalences between deformation problems.
\end{abstract}
\maketitle
\tableofcontents
\section{Introduction}
In this paper we develop the basic infinitesimal deformation theory of
\emph{abelian categories}.  This theory yields a natural
generalization of the well-known deformation theory of algebras
developed by Gerstenhaber \cite{GS,GerstI,GerstII}. In a subsequent paper
\cite{LVdBII} we will develop the corresponding obstruction theory in
terms of a suitable notion of Hochschild cohomology.

Deformation theory of abelian categories is important for
non-commutative algebraic geometry. One of the possible goals of
non-commutative algebraic geometry is to understand the
abelian (or triangulated) categories which have properties close to
those of the (derived) category of (quasi-)coherent sheaves on a
scheme. One is particularly interested in those
properties which are preserved under suitable deformations.  The
deformation theory of (abstract) triangulated categories seems
at this point somewhat elusive (due to the unclear status of the
currently accepted axioms) but, as we will show in this paper, there is a
perfectly good deformation theory for abelian categories.

As in any deformation theory we need some kind of flatness in order to
control the deformed objects.  Therefore the first contribution of
this paper is a notion of flatness for abelian
categories (see \S \ref{ref-3-31}). To the best of our knowledge this
definition is new.

In the rest of this introduction $R$ is a commutative coherent ring.
We will consider $R$-linear abelian categories. Informally these may
be viewed as non-commutative schemes over the (commutative) affine
base scheme $\Spec R$.

Let $\ccc$ be an $R$-linear abelian category. Our notion of flatness
has the following properties:
\begin{enumerate}
\item if $A$ is an $R$-algebra then $\Mod(A)$ is flat if and only if $A$
is flat over $R$;
\item $\ccc$ is flat if and only if $\ccc^{\op}$ is flat;
\item if $\ccc$ has enough injectives, $\ccc$ is flat  if and only
if injectives are flat~\cite{AZ2} in $\ccc^{\op}$;
\item if $\ccc$ is essentially small, $\ccc$ is flat if and only if
     $\Ind(\ccc)$ is flat (recall that $\Ind(\Cscr)$ is the formal
     closure of $\Cscr$ under filtered colimits (see \ref{ref-2.2-4}), it is a
     category with enough injectives);
\item flatness is stable under ``base change'' (see below).
\end{enumerate}
By enlarging the universe (see \S\ref{ref-2.1-0}) we may assume that
any category is small. Therefore, in principle, we could take properties (3)
and (4) as the definition of flatness. However this would make the
self duality (property (2)) very obscure. Our definition
of flatness is somewhat more technical (see \S\ref{ref-3-30}) but it is
manifestly left right symmetric.

Of fundamental importance in algebraic geometry is the concept of base
change. There is a natural substitute for this notion in our setting. Consider
a morphism $\theta:R\r S$ between coherent rings such that $S$ is finitely
presented over $R$ and let
$\Cscr$ be an $R$-linear category. The category $\Cscr_S$ is
the category of $S$-objects in $\ccc$, i.e. the pairs $(C,\varphi)$ where
$C\in \Ob(\Cscr)$ and $\varphi:S\r \Cscr(C,C)$ is an $R$-algebra map.
Intuitively $\Cscr_S$ the base extension of the ``non-commutative scheme''
$\Cscr$ to $\Spec S$.
We show in \S\ref{ref-4-41} that base change is compatible with various
natural constructions such as $\Ind$ and $\Mod$.

Assume now that $\theta$ is surjective such that
$I\overset{\text{def}}{=}\ker \theta$ satisfies $I^n=0$ for some $n$.
The surjectivity of $\theta$ implies that $\Cscr_S$ is a full abelian
subcategory of $\Cscr$.

Let $\Dscr$ be a flat $S$-linear category. A flat $R$-deformation of
$\Dscr$ is roughly speaking a flat lift of $\Dscr$ along the functor
$(-)_S$. In \S\ref{ref-6-64} we show that some of the basic properties of
abelian categories are preserved under flat deformation.
    More precisely, we show that the following
properties of  $\ddd$ lift to a flat  deformation.
\begin{enumerate}
\item $\ddd$ is essentially small;
\item $\ddd$ has enough injectives;
\item $\ddd$ is a Grothendieck category (i.e. a cocomplete abelian
category with a generator and exact filtered colimits);
\item $\ddd$ is a locally coherent Grothendieck category (i.e. $\ddd$
     is Grothendieck and is generated by a small abelian subcategory of
     finitely presented objects).
\end{enumerate}
In addition we show (see Theorem \ref{ref-8.5-119}) that up to equivalence
the number of flat deformations of an essentially small, respectively a
Grothendieck category is small.

Flatness is necessary for some of these properties. For example if
$k$ is a field then there are non-flat deformations of $\Mod(k)$
which do not have enough injectives (Example \ref{ref-6.17-83})
and furthermore the number of non-flat deformations is not small (Remark \ref{ref-8.6-120}).

In \S\ref{ref-7-104} we discuss the compatibility of localization with
deformations. Among other things we show that a deformation of $\Dscr$ gives
rise to deformations of all its localizations (Theorem \ref{ref-7.1-106}).

As a preparation to the sequel of this paper in
which we will develop the obstruction theory of abelian categories, we study
the associated deformation functors in \S\ref{ref-8-115}.  By
$\Def_\Dscr(R)$ we denote the flat $R$-deformations of $\Dscr$.

We have the following results.
\begin{enumerate}
\item If $\Dscr$ is an essentially small flat $S$-linear category then
     there is an equivalence between $\Def_{\Dscr}(R)$ and
     $\Def_{\Ind(\Dscr)}(R)$.
\item If $\Dscr$ is a locally coherent flat $S$-linear Grothendieck
     category then there is an equivalence between $\Def_{\Dscr}(R)$ and
     $\Def_{\Fp(\Dscr)}(R)$ where $\Fp(\Dscr)$ is the full (abelian) subcategory
of $\Dscr$ consisting of the finitely presented objects.
\end{enumerate}

In order to describe some more results
of  Section \S\ref{ref-8-115} we also need to introduce deformations of
general
$R$-linear categories. This is done by considering such categories as
``rings with several objects''~\cite{Mi}. We denote
by  $\ddef_{\frak{b}}(R)$   the flat $R$-deformations of $\frak{b}$ as
$R$-linear category.  Note that deforming an $R$-linear
abelian category $\Dscr$ \emph{as abelian category} is completely
different from deforming $\Dscr$ as $R$-linear category.

We prove the following results.
\begin{enumerate}
\setcounter{enumi}{2}
\item If $\frak{b}$ is an essentially small flat $S$-linear category
     then then there is an equivalence between $\ddef_{\frak{b}}(R)$ and
$\Def_{\Mod(\frak{b})}(R)$ where $\Mod(\frak{b})$ is the category of
covariant additive functors from $\frak{b}$ to $\Ab$. In particular
if $B$ is a flat $S$-algebra then there is an equivalence between
the deformations of $B$ and the deformations of $\Mod(B)$.
\item
If $\Dscr$ is a flat $S$-linear category with enough injectives then
there is an equivalence between $\Def_\Dscr(R)$ and $\ddef_{\Inj(\Dscr)}(R)$
where $\Inj(\Dscr)$ denotes the full (additive) subcategory of $\Dscr$
of injective objects.
\end{enumerate}
Property (3) shows that indeed our deformation theory generalizes
the deformation theory of algebras.

In  the final section of this paper we apply our methods
to the deformations of the category $\Mod(\Oscr_X)$ of sheaves of modules
over a ringed space
$(X,\Oscr_X)$.  For simplicity of exposition we assume here that $\Oscr_X$ is
a sheaf of
$k$-algebras where $k$ is a field. Assume that $X$ has a basis $\Bscr$
satisfying the following acyclicity condition:
\[
\forall U\in \Bscr: H^i(U,\Oscr_U)=0\qquad \text{for $i=1,2$}
\]

Let $\Oscr_\Bscr$ be the restriction of $\Oscr_X$ to $\Bscr$ and let
$\PMod(\Oscr_\Bscr)$ be the corresponding category of presheaves. 
We show that there is an equivalence
$$\Def_{\PMod(\Oscr_\Bscr)}(R) \cong \Def_{\Mod(\Oscr_X)}(R)$$

Let $\mathfrak{u}$ be the pre-additive category spanned by the (presheaf)
extensions by zero of the $\Oscr_U$ for $U\in\Bscr$. In other words we may
take $\Ob(\mathfrak{u})=\Bscr$ and we have
\[
\UUU(U,V) = \left\{
\begin{array}{r@{\hskip 1cm}l}
\ooo_X(U)& \mbox{if } U\subset V\\
0& \mbox{otherwise}
\end{array}
\right.
\]

Using property (3) above, it is easy to see that
$$\Def_{\PMod(\Oscr_\Bscr)}(R) \cong \ddef_{\UUU}(R)$$

These results confirm the fundamental insight of Gerstenhaber and
Schack \cite{GS1,GS2} that one should define the deformations of a
ringed space $(X,\Oscr_X)$ not as the deformations of $\Oscr_X$ as a
sheaf of $k$-algebras, but rather as the deformations of the
$k$-linear category $\frak{u}$ (or of the
``diagram''$(\bbb,\ooo_{\bbb})$ in case $X \in \bbb$).  These
``virtual'' deformations are nothing but the deformations of the
abelian category $\Mod(\Oscr_X)$.

\section{Preliminaries}
\subsection{Universes}
\label{ref-2.1-0}
It is well-known \cite{maclane} that category theory needs some extension of
the Zermelo-Fraenkel axioms of set theory (ZF).  One
possible extension is given by the G\"odel-Bernays axioms (GB) which
incorporates classes into set theory.  This makes it possible to introduce the
category
$\Set$ while at the same time avoiding Russel's paradox.

This solution is not entirely satisfying since for
example one would also like to talk about $\mathsf{Cls}$, the
category of all classes, and there is no room for this notion in GB.  In
particular, in the deformation theory of categories we consider below,
this seems to lead to foundational problems.

To solve such problems
Grothendieck introduced a more flexible extension of the
Zermelo-Fraenkel system: the theory of universes \cite{SGA4}.
The theory of universes  does not introduce new types of objects
but adds the universe axiom (U)
below.

A universe $\Uscr$ is a set with the following properties:
\begin{enumerate}
\item if $x\in \Uscr$ and if $y\in x$ then $y\in \Uscr$;
\item if $x,y\in\Uscr$ then $\{x,y\}\in \Uscr$;
\item if $x\in \Uscr$ then the powerset $\Pscr(x)$ of $x$ is in $\Uscr$;
\item if $(x_i)_{i\in I}$ is a family of objects in $\Uscr$ indexed
by an element of $\Uscr$ then $\bigcup_{i\in I} x_i\in \Uscr$;
\item if $U \in U$ and $f: U \lra \uuu$ is a function, then $\{ f(x) \,|\, x
\in U\} \in \uuu$.
\end{enumerate}
Note that $(x,y)$ is defined as $\{\{x,y\},x\}$ and hence if $x,y\in\Uscr$
then so is $(x,y)$.
A universe $\uuu$ with $\N \in \uuu$ is itself a model for
ZF.

As the only known non-empty universe only contains
finite sets we need the following new axiom:
\begin{enumerate}
\item[(U)] every set is the element of a universe.
\end{enumerate}

In this paper, we will work with ZFCU (the ZF axioms + the axiom of
choice + the universe axiom). By requiring $\{x,\N\} \in
\uuu$, every set $x$ is the element of a universe containing $x$ and
$\N$.  From now on, by a universe we will always mean a universe
containing $\N$. In particular, every such universe is itself a model
for ZFC.

We now  recall some terminology.
\begin{definition}
\begin{enumerate}
\item A set or cardinal is \emph{$\Uscr$-small} if it has the same
       cardinality as an element of $\Uscr$.
\item A category $\ccc$ consists of a set(!) of ``objects'' $\Ob(\ccc)$ and
a set of ``arrows'' $\Mor(\ccc)$ with the usual extra structure.
\item $\Uscr{-}\Set$ is the category whose objects consist of elements
of $\Uscr$ and whose $\Hom$-sets are just the standard $\Hom$'s between
sets. Likewise if $\Escr$ is a ``structure'' \cite{Bourbakistructure} (e.g.
abelian groups or rings) then $\Uscr{-}\Escr$ is the category of
$\Escr$-objects whose underlying set is in $\Uscr$. In particular,
$\uuu{-}\Cat$ (resp.
$\uuu{-}\Gd$) is the category of categories (resp. groupoids) $\ccc$ with
$\Ob(\ccc) \in
\uuu$ and $\Mor(\ccc) \in \uuu$ and the usual $\Hom$'s between categories.
\item A category is \emph{$\Uscr$-small} if both its objects and arrows are
       $\Uscr$-small sets.
\item A category is \emph{essentially $\Uscr$-small} if it equivalent to
a $\Uscr$-small category.
\item A category is a \emph{$\Uscr$-category} if it has $\Uscr$-small
$\Hom$-sets.
\item An abelian $\Uscr$-category $\Cscr$ is \emph{$\Uscr$-Grothendieck} if
       $\Cscr$ has a generator; $\Uscr$-small colimits exist in $\Cscr$ and
       $\Uscr$-small filtered colimits are exact.
\end{enumerate}
\end{definition}

\begin{remark}\label{ref-2.2-1}
If $E \in \Ob(\uuu{-}\eee)$, then $E \in \uuu$ since $E$ is described by an
element of $\uuu$. For example, $\ccc \in
\Ob(\uuu{-}\Cat)$ is described by an element of $\uuu^{\times 6} \subset \uuu$.
\end{remark}

\begin{lemma}\label{ref-2.3-2}
For $\ccc, \ddd \in \uuu{-}\Cat$, $\Hom(\ccc,\ddd) \in
\uuu{-}\Cat$.\hfill \qed
\end{lemma}

\begin{remark}\label{ref-2.4-3}
The axiom of choice allows us to replace a $\uuu$-category $\ccc$ with an
isomorphic category $\ccc'$ with $\Ob(\ccc') = \Ob(\ccc)$ and $\ccc'(C,D) \in
\uuu$ for every $C, D \in \Ob(\ccc)$.
In particular, if $\ccc$ is a (pre-additive) $\uuu$-category, we can define
representable functors $$\ccc(C,-): \ccc \lra \uuu{-}\Set (\uuu{-}\Ab):
D \longmapsto \ccc'(C,D)$$
where $\ccc'$ is as above.
\end{remark}

The universe axiom is the basis for the very useful
``extension of the universe'' principle.  I.e.\ by selecting a large enough
universe we may assume that any individual category is small.  The theory of
universes comes at a price however, namely the dependence of the
notations on the chosen universe.  Since this is rather tedious one
usually fixes the universe in advance and then drops it from the
notations except when invoking the extension of the universe principle. We
will follows these conventions in this paper.

\subsection{Some constructions depending on the universe}
\label{ref-2.2-4}\label{ref-2.2-5}
If $\mathfrak{a}$ is a pre-additive category then we denote by
$\Uscr{-}\Mod(\mathfrak{a})$ the category of covariant
additive functors
from $\mathfrak{a}$ to $\Uscr{-}\Ab$. It is
easy to see that if $\mathfrak{a}$ is essentially $\uuu$-small then
$\Uscr{-}\Mod(\mathfrak{a})$ is a $\Uscr$-Grothendieck category. If $\AAA$ is
a $\uuu$-category, $\uuu{-}\Mod(\AAA)$ contains functors $\AAA(A,-): \AAA \lra
\uuu-\Ab$ for
$A
\in
\AAA$ (see Remark \ref{ref-2.4-3}). In this case we define $\uuu-\mmod(\AAA)$ as
the full subcategory of
$\uuu{-}\Mod(\AAA)$ containing all functors that can be written as cokernels of
maps $\bigoplus_{i=1}^m
\mathfrak{a}(A_i,-)\r\bigoplus_{j=1}^n\mathfrak{a}(B_j,-)$. If $\uuu
\subset \vvv$ is an inclusion of universes, then Yoneda's lemma yields an
equivalence of categories
$\Uscr{-}\mod(\mathfrak{a})\lra\Vscr{-}\mod(\mathfrak{a})$.

If $\Cscr$ is an essentially small $\Uscr$-category then $\Uscr{-}\Ind(\Cscr)$
is
the full subcategory of $\Uscr{-}\Mod(\Cscr)$ consisting of left
exact functors.
It is well-known that $\Uscr{-}\Ind(\Cscr)$ is a $\Uscr$-Grothendieck category.
The objects in $\Uscr{-}\Ind(\Cscr)$ may be written as formal $\Uscr$-small
filtered colimits of objects in $\Cscr$ and the $\Hom$-sets are computed by the
rule
\[
\Hom_{\Uscr{-}\Ind(\Cscr)}(\colim_{i\in I} A_i,\colim_{j\in J} B_j)=
\lim_{i\in I}\colim_{j\in J} \Hom_{\Cscr}(A_i,B_j).
\]

An object $C$ in a $\uuu$-category $\ccc$ is $\uuu$-finitely presented
if the functor $\ccc(C,-): \ccc \lra \uuu{-}\Set$ (see Remark
\ref{ref-2.4-3}) preserves $\uuu$-small filtered colimits. We define
$\uuu{-}\Fp(\ccc)$ as the full subcategory of $\ccc$ containing
precisely the $\uuu$-finitely presented objects. It is well-known that
if $\ccc$ contains a $\uuu$-small full subcategory $\GGGG$ of
$\uuu$-finitely presented generators of $\ccc$, then
$\uuu{-}\Fp(\ccc)$ is the finite colimit closure of $\GGGG$ in $\ccc$
\cite{HCA2} and in particular is essentially $\uuu$-small.  If $\AAA$
is essentially $\uuu$-small, it is well-known and easy to see that
$\uuu{-}\Fp(\uuu{-}\Mod(\AAA)) = \uuu{-}\mmod(\AAA)$.

A $\Uscr$-Grothendieck category $\Cscr$ is locally coherent if it has a
$\uuu$-small set of $\uuu$-finitely presented generators and
$\uuu{-}\Fp(\Cscr)$ is a (necessarily essentially $\Uscr$-small) abelian
category.

If $\Cscr$ is a locally
coherent $\Uscr$-Grothendieck category then the natural functor
\[
\Uscr{-}\Ind(\uuu{-}\Fp(\Cscr))\lra\Cscr
\]
is an equivalence of categories.
If $\Cscr$ is essentially $\uuu$-small then the natural functor
\[
\Cscr\lra \uuu{-}\Fp(\Uscr{-}\Ind(\Cscr))
\]
is an equivalence as well.
\newtheorem*{convention}{Convention}
\begin{convention}
     From now on we work with a fixed universe $\Uscr$. All categories will
be $\Uscr$-categories. The notions of small and essentially small are
with respect to $\Uscr$. The same holds for the
notion of a Grothendieck category. The symbols $\Mod$, $\Ind$, $\Fp$ and
$\mmod$ are implicitly prefixed by $\Uscr$. Individual objects
such as abelian groups, rings, modules are, unless otherwise specified,
assumed to be $\Uscr$-small.
\end{convention}

\subsection{$R$-linear abelian categories}\label{ref-2.3-6}

Consider a commutative ring $R$. Recall that an \emph{$R$-linear
       category} is a pre-additive category $\AAA$ together with a ring map
$\rho: R \rightarrow \mathrm{Nat}(1_{\AAA},1_{\AAA})$. $\rho$ induces
a ring map $\rho_A: R \lra \AAA(A,A)$ for every object $A$ and an
action of $R$ on every $\Hom$-set. This leads to the equivalent
definition of an $R$-linear category as a category enriched in the
category \RMod of  $R$-modules. A pre-additive category is of
course the same as a $\Z$-linear category.

Unless otherwise stated, we will assume the ring $R$ to be coherent. Let
$\Rmod$ denote the full abelian subcategory of $\RMod$ of finitely presented
$R$-modules. Consider an abelian $R$-linear category $\ccc$. For every object
$C$ of $\ccc$ we obtain a (up to a canonical natural isomorphism) unique finite
colimit preserving functor
\begin{equation}(-
\otimes_R C): \Rmod
\lra \ccc
\end{equation}with $R \otimes _R C = C$.
We can construct its left derived functors $\Tor_i^R(-,C)$
using projective resolutions in $\mod(R)$.
In fact, we naturally have ($R$-bilinear) bifunctors
$$\Tor_i^R(-,-):
\Rmod
\times
\ccc
\lra \ccc.$$ The functors $(X \otimes_R -)$ are finite colimit preserving.
It is easily seen (using a fixed free resolution of $X$) that the functors
$\Tor_i^R(X,-)$ form a homological $\delta$-functor.

In a completely analogous way we define $$\Hom_R(-,C): (\Rmod)\op \lra \ccc$$
as the unique finite limit preserving functor with
$\Hom_R(R,C) = C$, and taking its right derived functors (again using
projective resolutions in
$\Rmod$) we obtain ($R$-bilinear) bifunctors $$\Ext_R^i(-,-):
(\Rmod)\op
\times
\ccc
\lra
\ccc.$$ The functors
$\Hom_R(X,-)$ are finite limit preserving. The functors $\Ext_R^i(X,-)$ form a
cohomological $\delta$-functor.

Note that $\Hom_R(X,-): \ccc \lra \ccc$ is in fact nothing but the opposite of
the tensor functor $(X \otimes_R -): \ccc\op \lra \ccc\op$ of $\ccc\op$.

\begin{definition}\label{ref-2.5-7}\cite{AZ2}
An object $C$ of $\ccc$ is \emph{flat (over $R$)} if the functor $(-
\otimes_R C):
\Rmod \lra \ccc$ is exact. Dually, $C$ is \emph{coflat (over $R$)} if
$\Hom_R(-,C):
\Rmod \lra \ccc$ is exact.
\end{definition}

\begin{remark}\label{ref-2.6-8} Flatness of $C$ is equivalent to the vanishing
of
$\Tor_1^R(-,C)$ and to the vanishing of all $\Tor_i^R(-,C)$ for $i \geq 1$.
Dually, coflatness of
$C$ is equivalent to the vanishing of $\Ext^1_R(-,C)$ and to the vanishing of
all
$\Ext^i_R(-,C)$ for $i \geq 1$. Consequently, in an exact sequence $0 \lra
A
\lra B \lra C \lra 0$, if $B$ and $C$ are flat, the same holds for $A$, and
dually, if $A$ and $B$ are coflat, the same holds for $C$.
\end{remark}

We will now list some useful facts concerning the $\Hom_R$ and $\otimes_R$
functors.

\begin{proposition}\label{ref-2.7-9}
Consider an $R$-linear functor $F: \ccc \lra \ddd$ between abelian
$R$-linear categories. If $F$ is exact, then for $X
\in
\Rmod$ and $C \in \ccc$
\begin{equation}\label{ref-2.2-10}
\Ext^i_R(X,F(C)) = F(\Ext^i_R(X,C)).
\end{equation}
If $F$ is only left exact, then \eqref{ref-2.2-10} holds for $i=0$.
\begin{proof}
This follows if we replace $X$ by a free resolution.
\end{proof}
\end{proposition}

\begin{corollary}\label{ref-2.8-11}
Consider an $R$-linear category $\ccc$ and a small category $\iii$. If
$\iii$-colimits are exact in $\ccc$, $\Ext^i_R(X,-)$ preserves them for $X \in
\Rmod$.
\begin{proof}
By taking $F$ to be the
functor $\colim:
\Fun(\iii,\ccc) \lra \ccc\hfill$, this follows from Proposition \ref{ref-2.7-9}.
\end{proof}
\end{corollary}

\begin{proposition}\label{ref-2.9-12}
For $X, Y \in \Rmod$ and $C, D, E \in \ccc$, the following hold:

\begin{enumerate}\setcounter{enumi}{-1}
\item \label{ref-0-13}\label{ref-0-14} $\ccc (X \otimes_R C,D) \cong \Hom_R (X,
\ccc(C,D))
\cong \ccc (C,
\Hom_R(X,D));$
\item \label{ref-1-15} $\Hom_R(X \otimes_R Y,C) \cong \Hom_R(X,\Hom_R(Y,C))$;
\item \label{ref-2-16} $(X \otimes_R Y)\otimes_R C \cong X \otimes_R (Y
\otimes_R C)$;
\item \label{ref-3-17}if $C$ is coflat, $X \otimes_R \Hom_R(Y,C) \cong
\Hom_R(\Hom_R(X,Y),C)$;
\item \label{ref-4-18}if $D$ is flat, $\Hom_R(X,Y \otimes_R D) \cong
\Hom_R(X,Y) \otimes_R D$;
\item \label{ref-5-19}if $E$ is injective, $X \otimes_R \ccc(C,E) \cong
\ccc(\Hom_R(X,C),E)$;
\item \label{ref-6-20}if $E$ is injective and $C$ is coflat, then $\ccc(C,E)$
is flat;
\item \label{ref-7-21}if $E$ is injective and $C$ is flat, then $\ccc(C,E)$ is
coflat;
\item \label{ref-8-22}if $\Ext^1(C,Z\otimes_R E) = 0$ for all
$Z$ in \Rmod  and
$E$ is flat, then $X \otimes_R \ccc(C,E) = \ccc(C,X \otimes_R E)$ and
$\ccc(C,E)$ is flat.
\hfill\qed
\end{enumerate}
\end{proposition}

Consider a finitely generated ideal $I$ of $R$. For every $C \in \ccc$, we
obtain a map $I \otimes_R C \lra C$, which is not necessarily mono unless
$C$ is flat. We will denote the image of this map by $IC$. Dually, we will
denote the image of the map $C \lra \Hom_R(I,C)$ by $CI$. It is easy to see
that $(IJ)C = I(JC)$ as subobjects of $C$.

\medskip

If $(\AAA, \rho)$ is an $R$-linear category and $F \in \Mod(\mathfrak{a})$,
then the abelian groups $F(A)$ inherit an $R$-module structure from the maps
$\rho(r): A \lra A$. In fact, $\Mod(\mathfrak{a})$ is isomorphic to the
$R$-linear category of $R$-linear functors from $\AAA$ into $\RMod$ (see
also \S\ref{ref-4-40}).

\begin{proposition}\label{ref-2.10-23}
     For an essentially small $R$-linear category $\AAA$, the $\Hom_R$
     and $\otimes_R$ functors on $\Mod(\AAA)$ are computed pointwise.
     I.e.
\begin{align*}
(X\otimes_R F)(A)&=X\otimes_R F(A)\\
\Hom_R(X,F)(A)&=\Hom_R(X,F(A))
\end{align*}
for $X\in \mod(R)$, $F\in \Mod(\AAA)$ and $A\in \AAA$.
In particular, $F$ in $\Mod(\AAA)$ is flat (resp. coflat) if and only if all
its values $F(A)$ are flat (resp. coflat) in $\Mod(R)$.\hfill \qed
\end{proposition}

\subsection{Derived functors and ind-objects}

Some arguments in this paper are based on extending derived functors to
ind-objects (see \S \ref{ref-2.2-4}). Therefore in this section we discuss
some of the relevant properties in this regard.

\begin{definition}\label{ref-2.11-24}
\begin{enumerate}
\item A functor $F: \aaa \lra \bbb$ between an arbitrary category $\aaa$ and a
       pre-additive category $\bbb$ is called \emph{effaceable}\cite{Groth1} if
for every $A
       \in \aaa$ there is a monomorphism $u: A \lra A'$ with $F(u) = 0$.
       $F$ is called \emph{co-effaceable} if for every $A$ there is an
       epimorphism $u: A' \lra A$ with $F(u) = 0$.
\item A functor $F: \aaa \lra \Ab$ is called \emph{weakly
effaceable} if
for each $A \in \aaa$ and $x \in F(A)$ there is a monomorphism $u: A \lra A'$
with $F(u)(x) = 0$.
\end{enumerate}
\end{definition}

\begin{proposition}\label{ref-2.12-25}\cite{Groth1}
Let $F: \aaa \lra \bbb$ be an additive functor between abelian categories. If
$\aaa$ has enough injectives then the following are equivalent:

\begin{enumerate}
\item $F$ is effaceable;
\item $F(I) = 0$ for all injectives $I$ of $\aaa$.\hfill\qed
\end{enumerate}
\end{proposition}

Consider an extension to ind-objects
$$\xymatrix{{\ccc} \ar[r]\ar[d]_F & {\Ind(\ccc)}
\ar[d]^{\Ind(F)}
\\
\ddd \ar[r]& {\Ind(\ddd)}}$$
where $\Cscr$, $\Dscr$ are essentially small abelian categories.
\begin{proposition}\label{ref-2.13-26}
The following are equivalent:
\begin{enumerate}
\item $F$ is effaceable;
\item $\Ind(F)$ is effaceable.
\end{enumerate}
\begin{proof}
Suppose $\Ind(F)$ is effaceable and consider $C$ in $\ccc$. There is a
monomorphism $u: C \lra \colim_iC_i$ in $\Ind(\ccc)$ with $\Ind(F)(u) = 0$ for
certain $C_i$ in $\ccc$. Consider the
maps $s_j: C_j \lra \colim_iC_i$. Since $C$ is finitely presented, $u$
factorizes as $u = s_j \circ f$ for some $f: C \lra C_j$. Now $s_j$ is mapped
by $\Ind(F)$ onto the canonical map $F(C_j) \lra \colim_iF(C_i)$, and its
composition with $F(f)$ is zero. Since $F(C)$ is finitely presented in
$\Ind(\ddd)$, there is a map $s_{jk}:C_j \lra C_k$ such that $F(s_{jk}) \circ
F(f) = 0$, so $s_{jk} \circ f$ is a monomorphism that serves our purpose.

Now suppose $F$ is effaceable. By Proposition \ref{ref-2.12-25} it
suffices to prove that
$\Ind(F)(I) = 0$ for every injective $I = \colim_iC_i$ of $\Ind(\ccc)$. For
every $C_i$, we take a monomorphism $u_i: C_i \lra C_i'$ in $\ccc$ with $F(u_i)
= 0$. Since $I$ is injective, the maps $s_j:C_j \lra I$ factorize as $s_j = f_j
\circ u_j$. It follows that all the maps $\Ind(F)(s_j): F(C_j) \lra
\colim_iF(C_i)$ are zero and hence $\Ind(F)(I) = \colim_iF(C_i) = 0$.
\end{proof}
\end{proposition}

If $\ccc$ is an abelian category then $\Ext^i_{\ccc}(-,-)$ denotes the
\emph{Yoneda} $\Ext$-groups between $\ccc$-objects. If $\ccc$ has enough
injectives, then $\Ext^i_{\ccc}(C,-)_i$ are the derived functors of
$\ccc(C,-)$. The following proposition is presumably well-known, but we have
been unable to find a reference.

\begin{proposition}\label{ref-2.14-27} Assume that $\Cscr$ is essentially small.
For $A, B \in \ccc$, we have
\begin{equation} \label{ref-2.3-28}
\Ext^i_{\ccc}(A,B) = \Ext^i_{\Ind(\ccc)}(A,B).
\end{equation}
\begin{proof}
The formula \eqref{ref-2.3-28} is obviously true for $i = 0$. Then by
\cite[Lemma II.2.1.3]{Illusie}, it is sufficient to show that the functor
$\Ext^i_{\Ind(\ccc)}(-,B)$ restricted to $\ccc$ is weakly effaceable. This
follows from Proposition
\ref{ref-2.15-29} below with $G = \Ext^i_{\Ind(\ccc)}(-,B)$.
\end{proof}
\end{proposition}

\begin{proposition}\label{ref-2.15-29} Assume that $\Cscr$ is essentially small.
If in the following commutative diagram
$$\xymatrix{{\ccc\op} \ar[r]\ar[rd]_F & {(\Ind(\ccc))\op} \ar[d]^{G} \\
& {\Ab}}$$
$G$ is weakly effaceable, the same holds for $F$.
\begin{proof}
Consider $C$ in $\ccc$ and $x \in F(C)$. There is an epimorphism $f:
\colim_iC_i
\lra C$ in
$\Ind(\ccc)$ with $G(f)(x) = 0$ for certain $C_i$ in $\ccc$. Consider the
maps $s_j: C_j \lra \colim_iC_i$. We have that $C = \sum_i \Beeld(f \circ
s_i)$,
but since $C$ is finitely generated in $\Ind(\ccc)$, we find an epimorphism $f
\circ s_j: C_j \lra C$ in $\Ind(\ccc)$ (and thus in $\ccc$) with $F(f \circ
s_j)(x) = G(f \circ s_j)(x) = G(s_j)(G(f)(x)) = 0$.
\end{proof}
\end{proposition}

\section{Flatness for $R$-linear abelian
categories}\label{ref-3-30}\label{ref-3-31}

An $R$-linear category is called \emph{flat} if its $\Hom$-sets are flat in
$\Mod(R)$. In this section, we introduce a different notion of flatness for an
$R$-linear \emph{abelian} category $\ccc$. This notion has the following
properties:
\begin{enumerate}
\item $\ccc$ is flat if and only if $\ccc^{\op}$ is flat;
\item if $\ccc$ has enough injectives, $\ccc$ is flat if and only if
injectives are coflat;
\item if $\ccc$ is essentially small, $\ccc$ is flat if and only if
$\Ind(\ccc)$ is flat;
\item if $\AAA$ is an essentially small $R$-linear category, $\AAA$ is flat
(as an $R$-linear category) if and only if $\Mod(\AAA)$ is flat (as an
abelian $R$-linear category).
\end{enumerate}
By enlarging the universe (see \S\ref{ref-2.1-0}) we may assume that any
category is small. Therefore, in principle, we could take properties (3)
and (4) as the definition of flatness. However this would make the
self duality (property (2)) very obscure. Therefore we give a somewhat more
technical, but manifestly left right symmetric definition below.

Let $\ccc$ be an $R$-linear abelian category. The following result
is crucial for what follows.

\begin{proposition}\label{ref-3.1-32}
For $X \in \Rmod$, the following are equivalent:

\begin{enumerate}
\item $\Ext_R^1(X,-): \ccc \lra \ccc$ is effaceable;
\item $\Ext^1_R(X, \ccc(-,-)): \ccc\op \times \ccc \lra \Ab$ is weakly
effaceable;
\item $\Tor_1^R(X,-): \ccc \lra \ccc$ is co-effaceable.
\end{enumerate}

\begin{proof}
Since (2) is clearly self-dual, and (3) is the dual statement of (1), it
suffices to prove the equivalence of (1) and (2).
Take a projective resolution $P_{\cdot}$ of $X$ and let $Y$ denote the
kernel of $P_0 \lra X$. Note that since $$\Ext^1(X,\ccc(A,B)) = H_{1}\ccc(A,
\Hom_R(P_{\cdot},B)),$$ an element in this group can be represented
by a $\ccc$-map $A
\lra \Hom_R(P_1,B)$ such that the composition with $\Hom_R(P_1,B) \lra
\Hom_R(P_2,B)$ is zero.
Suppose (2) holds. Consider an arbitrary $\ccc$-map $A \lra
\Ext_R^1(X,B)$. We will show that there is a mono $B \lra B'$ such that $A \lra
\Ext^1_R(X,B)
\lra \Ext^1_R(X,B')$ is zero, which proves (1). Consider the pullback
$$\xymatrix{{A'} \ar[r]\ar[d] & A \ar[d] \\ {\Hom_R(Y,B)} \ar[r]&
{\Ext^1_R(X,B).}}$$
The composition $A' \lra \Hom_R(Y,B) \lra \Hom_R(P_1,B)$ represents an element
of
$\Ext^1(X,\ccc(A,B))$, hence by (2) we can take $A'' \lra A'$ and $B
\lra B'$ such that $A'' \lra A' \lra \Hom_R(Y,B) \lra \Hom_R(Y,B')$
factors over $\Hom_R(P_0,B') \lra \Hom_R(Y,B')$. It follows that the
composition $A'' \lra A' \lra \Hom_R(Y,B) \lra \Hom_R(Y,B') \lra
\Ext^1_R(X,B')$ is zero, and consequently $A \lra \Ext^1_R(X,B)
\lra \Ext^1_R(X,B')$ is zero too.

Now suppose (1) holds and consider an element in $\Ext^1(X,\ccc(A,B))$
represented by $A \lra \Hom_R(P_1,B)$. Take a monomorphism $B \lra B'$ such
that $\Ext^1_R(X,B) \lra \Ext^1_R(X,B')$ is zero. Let $K$ denote the kernel
of $\Hom_R(P_1,B) \lra \Hom_R(P_2,B)$ and let $L'$ denote the image of
$\Hom_R(P_0,B') \lra \Hom_R(P_1,B')$. We obtain the following
factorization
$$\xymatrix{ {A'} \ar[dd] \ar[r]& A \ar[d] \ar[dr] \\
& K \ar[r]\ar[d] & {\Hom_R(P_1,B)} \ar[d] \\ {\Hom_R(P_0,B')} \ar[r]& {L'}
       \ar[r]& {\Hom_R(P_1,B'),}}$$
in which $A'$ is defined as a pullback. Clearly, the composition $A' \lra
\Hom_R(P_1,B')$ represents zero in $\Ext^1(X,\ccc(A',B'))$.
\end{proof}
\end{proposition}

\begin{definition}\label{ref-3.2-33}
The category $\ccc$ is called \emph{flat (over $R$)} if every $X \in \Rmod$
satisfies the equivalent conditions of Proposition \ref{ref-3.1-32}.
\end{definition}
Note that the notion of flatness is independent of the choice of universe.

\begin{proposition}\label{ref-3.3-34}
The following are equivalent:
\begin{enumerate}
\item $\ccc$ is flat over $R$;
\item $\ccc\op$ is flat over $R$.
\end{enumerate}
\begin{proof}
Immediate from Proposition \ref{ref-3.1-32}.
\end{proof}
\end{proposition}

\begin{proposition}\label{ref-3.4-35}
If $\ccc$ has enough injectives, the following are equivalent:
\begin{enumerate}
\item $\ccc$
is flat over $R$;
\item injectives in $\ccc$ are coflat.
\end{enumerate}
\begin{proof}
This follows from Proposition \ref{ref-2.12-25} with $F = \Ext^1_R(X,-)$.
\end{proof}
\end{proposition}

\begin{proposition}\label{ref-3.5-36}\label{ref-3.5-37}
If $\ccc$ is flat (over $R$), the following hold:
\begin{enumerate}
\item For every $X \in \mmod(R)$ and $i \geq 1$, $\Ext_R^i(X,-)$ is
effaceable;
\item For every $X \in \mmod(R)$, the cohomological $\delta$-functor
$(\Ext_R^i(X,-))_{i\ge 0}$ is universal \cite{Groth1};
\item for every $X \in \mmod(R)$ and $i \geq 1$, $\Tor^R_i(X,-)$ is
co-effaceable.
\item for every $X \in \mmod(R)$, the homological $\delta$-functor
$(\Tor_i^R(X,-))_{i \ge 0}$ is universal.
\end{enumerate}
\begin{proof}
(1) and (3) follow by dimension shifting in the first argument. (2) and (4)
follow from (1) and (3) by \cite{Groth1}.
\end{proof}
\end{proposition}

The category of ind-objects $\Ind(\ccc)$ (see \S\ref{ref-2.2-4}) is obviously
$R$-linear. It follows from Proposition
\ref{ref-2.7-9} that the functor $\ccc \lra \Ind(\ccc)$ commutes with
$\Tor^R$ and $\Ext_R$. So $\Tor_i^{\Ind(\ccc),R}(X,-)$ and
$\Ext^i_{\Ind(\ccc),R}(X,-)$ are the extensions to
ind-objects of $\Tor_i^{\ccc,R}(X,-)$ and $\Ext^i_{\ccc,R}(X,-)$.

\begin{proposition}\label{ref-3.6-38}
Assume that $\Cscr$ is essentially small. The following are equivalent:
\begin{enumerate}
\item $\ccc$ is flat over $R$;
\item $\Ind(\ccc)$ is flat over $R$.
\end{enumerate}
\begin{proof}
This follows immediately from Proposition \ref{ref-2.13-26} with $F =
\Ext^1_{\ccc,R}(X,-)$ and
$\Ind(F) = \Ext^1_{\Ind(\ccc),R}(X,-)$.
\end{proof}
\end{proposition}

\begin{proposition}\label{ref-3.7-39}
       For an essentially small $R$-linear category $\AAA$, the following
       are equivalent:
\begin{enumerate}
\item $\AAA$ is flat over $R$;
\item the abelian category $\AAAMod$ is flat over $R$.
\end{enumerate}
\begin{proof}
First, note that for every $E$ in $\Mod(\AAA)$ and for every $A$ in $\AAA$,
$$\begin{aligned}
\Mod(R)(X, E(A)) &= \Mod(R)(X,
\Mod(\AAA)(\AAA(A,-),E))\\ &= \Mod(\AAA)(X \otimes_R \AAA(A,-),E).
\end{aligned}$$
By Proposition \ref{ref-2.10-23}, $\Mod(\AAA)$ is flat if and only
if for every
injective $E$ in $\Mod(\AAA)$ and every $A$ in $\AAA$, $E(A)$ is coflat. Also
by Proposition \ref{ref-2.10-23}, $\AAA$ is flat if and only if every functor
$\AAA(A,-)$ for $A$ in $\AAA$ is flat. And since $\Mod(\AAA)$ has enough
injectives, this is equivalent to requiring that for every injective $E$ in
$\Mod(\AAA)$ the functor $\Mod(\AAA)(-,E) \circ (- \otimes_R \AAA(A,-)):
\mmod(R) \lra \Mod(R)$ is exact.
The statement then follows from the computation above.
\end{proof}
\end{proposition}

\section{Base change}
\label{ref-4-40}\label{ref-4-41}
In this section we study a natural substitute for the notion of base change in
algebraic geometry, and we show that it is compatible with constructions
such as $\Ind$ and $\Mod$.

We fix a homomorphism of commutative rings
$\theta: R
\lra S$. For
$M
\in
\Mod(S)$,
$\ovl{M}$ denotes $M$ considered as an $R$-module using $\theta$. For an
$S$-linear category $\BBB$,
$\ovl{\mathfrak{b}}$ is the $R$-linear category with $\Ob(\ovl{\mathfrak{b}})=
\Ob({\mathfrak{b}})$ and $\ovl{\mathfrak{b}}(B,B')=\ovl{\mathfrak{b}(B,B')}$.
More often than not we will  just write $\mathfrak{b}$ for
$\ovl{\mathfrak{b}}$.

For an
$R$-linear category $\AAA$, $S \otimes_R\AAA$ is the $S$-linear category with
$\Ob(S \otimes_R \AAA) = \Ob(\AAA)$ and $(S \otimes_R \AAA)(A,A') = S
\otimes_R
\AAA(A,A')$.

The functor $S \otimes_R (-)$ is left adjoint to  $\ovl{(-)}$ in the sense that
for an
$R$-linear category $\AAA$ and an
$S$-linear category $\BBB$, there is an isomorphism, natural in
$\mathfrak{a}$, $\mathfrak{b}$:
\[
\overline{\Add(S)(S \otimes_R\AAA,\BBB)}\cong\Add(T)(\AAA,\ovl{\BBB})
\]
as $R$-linear categories where $\Add(T)$ denotes the $T$-linear functors.

For an $R$-linear category $(\ccc, \rho)$, let $\ccc _S$ be the following
$S$-linear category of \emph{$S$-objects}: the objects of $\ccc _S$ are couples
$(C,\varphi)$ where $C$ is an object of $\ccc$ and $\varphi :S \lra
\ccc(C,C)$ is a ring map with $\varphi \circ \theta = \rho_C$. The
morphisms of $\ccc _S$ are the obvious compatible $\ccc$-morphisms.
Clearly $\ccc _S$ becomes $S$-linear using the ring maps $\varphi$.
The  functor $(-)_S$  is
right
adjoint to $\overline{(-)}$ in the sense that for an $R$-linear category $\ccc$
and an
$S$-linear category $\ddd$,
there is an isomorphism, natural in $\Cscr$, $\Dscr$:
$$\ovl{\Add(S)(\ddd,\ccc_S)} \cong \Add(R)(\ovl{\ddd},\ccc)$$
of $R$-linear categories.

\begin{remark}
It is readily seen that $(\Mod(R))_S \cong \Mod(S)$, hence we
deduce the  isomorphism $\Add(S)(\BBB,\Mod(S)) \cong
\Mod(\BBB)$ for an $S$-linear category $\BBB$ which was already mentioned
just before Proposition \ref{ref-2.10-23}.
\end{remark}

\begin{proposition}\label{ref-4.2-42}
     The forgetful functor $F: \ccc_S \lra \ccc$ preserves, reflects and
     creates \cite{AHS} all limits and colimits.
\begin{proof}
We give the construction of limits in $\ccc_S$ that exist in $\ccc$. Let $G:
\iii \lra \ccc_S$ with $G(i) = (C_i,\varphi_i)$ be an arbitrary functor and let
$\lim_iC_i$ denote the limit of $F \circ G$. For every $s \in S$,
the maps $\varphi_i(s): C_i \lra C_i$ define a natural transformation $F
\circ G \lra F \circ G$ mapped by the limit functor to $\lim_i(
\varphi_i(s))$. It is easily seen that $(\lim_iC_i, \lim_i\varphi_i)$ where
$(\lim_i\varphi_i)(s) = \lim_i(\varphi_i(s))$ is a limit of $G$.
\end{proof}
\end{proposition}
When
$\ccc$ is abelian, the same holds for $\ccc_S$, and in this case
$\ccc_S \lra \ccc$ is exact. From now on we will only use the functor $(-)_S$
for abelian categories.

Assume now that $R\lra S$ is a ring morphism between coherent commutative
rings such that $S$ is finitely presented as $R$-module.
The functors $(S \otimes_R -): \ccc \lra \ccc$ and $\Hom_R(S,-): \ccc \lra
\ccc$ of Section \ref{ref-2.3-6} yield functors
$$(S
\otimes_R -): \ccc \lra \ccc_S$$ and $$\Hom_R(S,-): \ccc \lra \ccc_S.$$

It is easily seen that
\begin{proposition}\label{ref-4.3-43}\label{ref-4.3-44}
The functors
$(S \otimes_R -)$ and $\Hom_R(S,-)$ are respectively left and right adjoint to
$\ccc_S \lra \ccc$. In particular, $(S \otimes_R -)$ preserves projectives
and $\Hom_R(S,-)$ preserves injectives. \hfill\qed
\end{proposition}
The next result gives a connection between the functors $S \otimes_R(-)$ and
$(-)_S$.
\begin{proposition}
\label{ref-4.4-45}
Let $\mathfrak{a}$ be an essentially small $R$-linear category. Then
\begin{enumerate}
\item there is a commutative diagram
\[
\begin{CD}
\Mod(S \otimes_R\mathfrak{a})@>\cong>> \Mod(\mathfrak{a})_S\\
@V\alpha VV  @VV\beta V\\
\Mod(\mathfrak{a}) @= \Mod(\mathfrak{a})
\end{CD}
\]
where $\alpha$ is dual to $\mathfrak{a}\r S \otimes_R\mathfrak{a}$ and $\beta$
is the forgetful functor;
\item the left and right adjoints to the forgetful functor
$\Mod(\mathfrak{a})_S
\r \Mod(\mathfrak{a})$
are computed pointwise, i.e.
\begin{align*}
(S\otimes_R F)(A)&=S\otimes_R FA\\
\Hom_R(S,F)(A)&=\Hom_R(S,F(A))
\end{align*}
In particular $S\otimes_R \mathfrak{a}(A,-)=(S
\otimes_R \mathfrak{a})(A,-)$.\hfill\qed
\end{enumerate}
\end{proposition}

We have the following relation between $\Ind(-)$ and $(-)_S$:

\begin{proposition}
\label{ref-4.5-46}
Assume that $\Cscr$ is essentially small.
The obvious functor
\[
H:\Ind(\Cscr_S)\r (\Ind(\Cscr))_S
\]
is an equivalence.
\end{proposition}
\begin{proof}
Consider the following commutative diagram with obvious maps.
$$\xymatrix{\ccc \ar[r]& \Ind(\ccc) &\\
\ccc_S \ar[u]^F \ar[r] & \Ind{(\ccc_S)} \ar[u]^{\Ind(F)} \ar[r]_-{H} &
{(\Ind(\ccc))_S}
\ar[ul]}$$

Since $F$ is
faithful, the same holds for $\Ind(F)$ and thus for $H$.

To prove that $H$ is
full, consider a map $(f_i)_i: \colim_i(C_i,\varphi_i) \lra
\colim_j(D_j,\psi_j)$ with $f_i: C_i \lra D_{j_i}$ maps in $\ccc$ such that
$(\psi_j)_j \circ (f_i)_i = (f_i)_i \circ (\varphi_i)_i$ in $\Ind(\ccc)$. This
means that for every $i$ we obtain a diagram
$$\xymatrix{C_i \ar[r]^-{f_i} \ar[d]_{\varphi_i} & D_{j_i}
\ar[r]^-{g_i} \ar[d]^{\psi_{j_i}} & D_{k_i} \ar[d]^{\psi_{k_i}} \\
C_i \ar[r]_-{f_i} & D_{j_i} \ar[r]_-{g_i} & D_{k_i}}$$
in which the right square commutes and composition with $g_i$ makes the left
square commute. It follows that the maps $g_i \circ f_i$ belong to
$\ccc_S$ and
$(f_i)_i = (g_i \circ f_i)_i$.

To prove that $H$ is essentially surjective, consider an object
$C$ of $(\Ind(\ccc))_S$ (We omit the $S$-action $\varphi$ in our notation).
$\varphi$ induces an $(\Ind(\ccc))_S$-epimorphism $S \otimes_R C \lra C$ with
kernel $K$, for which we obtain another $(\Ind(\ccc))_S$-epimorphism $S
\otimes_R K \lra K$. It follows that $C$ is isomorphic to the cokernel of $S
\otimes_R K \lra S \otimes_R C$. Writing $C = \colim_i C_i$ as a filtered
colimit, we see that $S \otimes_R C = \colim_i (S \otimes_R C_i)$ belongs to
$\Ind(\ccc_S)$ and the same holds for $S \otimes_R K$. Since $\Ind(F)$ is
exact,  the cokernel of the map belongs to $\Ind(\ccc_S)$, as we wanted.
\end{proof}

Next we consider finitely presented objects.

\begin{proposition}\label{ref-4.6-47}\label{ref-4.6-48}\label{ref-4.6-49}
Assume that in $\ccc$ small filtered colimits are exact. If $(C,\varphi)$ is
finitely presented (resp finitely generated) in $\ccc_S$, then $C$ is finitely
presented (resp. finitely generated) in
$\ccc$. The obvious functor
$$H: \Fp(\ccc_S) \lra (\Fp(\ccc))_S$$ is an isomorphism.
\begin{proof}
Consider a finitely presented $(C,\varphi)$ in $\ccc_S$ and a filtered colimit
$\colim_iC_i$ in
$\ccc$. Making use of Corollary \ref{ref-2.8-11} and the fact that $\ccc_S
\lra \ccc$ reflects colimits, we may compute
$$\begin{aligned}
\ccc(C,\colim_iC_i) & = \ccc_S((C,\varphi), \Hom_R(S,\colim_iC_i)) \\
& = \ccc_S((C,\varphi),\colim_i\Hom_R(S,C_i)) \\
& = \colim_i\ccc_S((C,\varphi),\Hom_R(S,C_i)) \\
& = \colim_i\ccc(C,C_i),
\end{aligned}$$
hence $C$ is finitely presented in $\ccc$.
Now consider the following commutative diagram with obvious maps.
$$\xymatrix{\ccc \ar[r]& \Fp(\ccc) &\\
\ccc_S \ar[u]^F \ar[r] & \Fp{(\ccc_S)} \ar[u]^{\Fp(F)} \ar[r]_-{H} &
{(\Fp(\ccc))_S}
\ar[ul]}$$
$H$ is readily seen to be fully faithful and injective on objects, so it
remains to show that it is surjective on objects. Consider $(C,\varphi)$ in
$(\Fp(\ccc))_S$ and a filtered colimit $\colim_i(C_i,\varphi_i) = (\colim_i
C_i, \colim_i \varphi_i)$ in $\ccc_S$. The canonical map
$$\colim_i\ccc_S((C,\varphi),(C_i,\varphi_i)) \lra \ccc_S((C,\varphi),
(\colim_iC_i,\colim_i\varphi_i))$$ is obviously injective. We finish the proof
by showing that it is surjective. If
$f: C
\lra
\colim_i C_i$ defines a map in the codomain, $f$ factors over some $f_i: C \lra
C_i$ in $\ccc$. For each
$s \in S$, there is a diagram
$$\xymatrix{C \ar[r] \ar[d]_{\varphi(s)}& {C_i} \ar[r]
\ar[d]^{\varphi_i(s)} & {\colim_i C_i} \ar[d]^{\colim_i\varphi_i(s)} \\
C \ar[r] & {C_i} \ar[r]& {\colim_iC_i}}$$
in which the right hand square is not necessarily commutative. Selecting
finitely many generators $s_k$ of $S$ over $R$, and using the fact that $C$
is finitely presented in $\ccc$ and that the colimit is filtered, we can find
$C_i \lra C_j$ such that the right hand square with $C_j$ instead of $C_i$
commutes for every generator $s_k$ and hence also for every $s \in S$.
\end{proof}
\end{proposition}

We mention the familiar change of rings spectral sequences.
\begin{proposition}\label{ref-4.7-50}
     Let $X \in \Smod$, $C \in \ccc$, $A\in \Cscr_S$. Assume that $\Cscr$
     is flat.  There are first quadrant spectral sequences
\begin{gather}
E^{pq}_2:\Ext^p_{\Cscr_S}(A,\Ext^q_R(S,C))\Rightarrow \Ext^n_{\Cscr}(A,C)
\label{ref-4.1-51}\\
E^{pq}_2:\Ext^p_S(X,\Ext^q_R(S,C))\Rightarrow \Ext^n_R(X,C)\label{ref-4.2-52}
\end{gather}
\end{proposition}
\begin{proof}  By enlarging our universe we may assume that $\Cscr$ is small.
     Then we replace $\Cscr$ by $\Ind(\Cscr)$ in order to have a category
     with enough injectives. By Proposition \ref{ref-3.6-38} this enlarged
     $\ccc$ is still flat. Furthermore by Proposition \ref{ref-3.5-36}
$\Ext^q_R(S,-)$ is the derived functor of $\Hom_R(S,-)$.

We have
\begin{align}
       \Hom_{\Cscr_S}(A,\Hom_R(S,C))&=\Hom_{\Cscr}(A,C)
       \label{ref-4.3-53}\\
       \Hom_S(X,\Hom_R(S,C))&= \Hom_R(X,C) \label{ref-4.4-54}
\end{align}
Now $\Hom_R(S,-)$ preserves injectives and \eqref{ref-4.4-54} shows
that $\Hom_R(S,-)$ also preserves coflat objects. So $\Hom_R(S,-)$
sends injectives to acyclic objects for
$\Hom_{\Cscr_S}(A,-)$ and $\Hom_S(X,-)$.

Thus
\eqref{ref-4.1-51}\eqref{ref-4.2-52} are just the Grothendieck spectral
sequences \cite{Groth1} associated to \eqref{ref-4.3-53}
\eqref{ref-4.4-54}.
\end{proof}
Here are some properties that are preserved under base change.
\begin{proposition}\label{ref-4.8-55}
\label{ref-4.8-56}\label{ref-4.8-57}
\begin{enumerate}
\item If $\ccc$ has enough injectives, the same holds for $\ccc_S$;
\item if $(G_i)_i$ is a set of generators of $\ccc$, $(S \otimes_R G_i)_i$
is a set of generators of $\ccc_S$;
\item if $\ccc$ is Grothendieck, the same holds for $\ccc_S$;
\item  if
$\ccc$ is a locally coherent Grothendieck category, the same holds for
$\ccc_S$;
\item
if $\ccc$ is flat over $R$, then $\ccc_S$ is flat over $S$.
\end{enumerate}
\end{proposition}
\begin{proof}
Only (5) is not entirely clear from the above discussion.
Consider $C \in \ccc_S$ and $X \in \mmod(S)$. Take a $\ccc$-monomorphism $m: C
\lra C'$ with $\Ext^1_R(X,m) = 0$. Then $\Hom_R(S,m): \Hom_R(S,C)
\lra \Hom_R(S,C')$ is a
monomorphism. The spectral sequence \eqref{ref-4.2-52} yields
monomorphism $\Ext^1_S(X,\Hom_R(S,C))\longrightarrow \Ext^1_S(X,C)$
natural in $C$ and
hence $\Ext^1_S(X, \Hom_R(S,m)) = 0$.
The composition of monomorphims $C\r \Hom_R(S,C) \lra \Hom_R(S,C')$ is
the required effacing of $\Ext^1_S(X, C)$ in $\Cscr_S$.
\end{proof}
\section{Deformations of abelian categories}
In this section we introduce deformations of pre-additive and abelian
categories.  Basically these are lifts along the functors $S \otimes_R(-)$ and
$(-)_S$ for a homomorphism of commutative rings $\theta: R\r S$.  In
the sequel we will assume that $\theta$ is surjective with nilpotent
kernel but this is not necessary for the basic definitions.

We start with the pre-additive case.
\begin{definition}\label{ref-5.1-58}
For an $S$-linear category $\BBB$, an \emph{$R$-deformation} of
$\BBB$ is an $R$-linear category $\AAA$ together
with an
$R$-linear functor $\AAA \lra \overline{\BBB}$ that induces an equivalence
$S \otimes_R\AAA \lra \BBB$. A deformation $\AAA \lra \overline{\BBB}$ will
often
       be denoted by $\AAA \lra \BBB$ or simply by $\AAA$. If
$S \otimes_R\AAA
\lra \BBB$ is an isomorphism, the deformation is called \emph{strict}.
If $\AAA$ is flat over $R$, the deformation is called \emph{flat}.
\end{definition}
Now we consider the abelian case.
\begin{definition}\label{ref-5.2-59}
For an \emph{abelian} $S$-linear category $\ddd$, an \emph{$R$-deformation} of
$\ddd$ is an \emph{abelian} $R$-linear category $\ccc$
together with an $R$-linear functor $\overline{\ddd} \lra \ccc$ that induces an
equivalence
$\ddd \lra \ccc _S$. A deformation $\overline{\ddd} \lra \ccc$ will often be
denoted by $\ddd \lra \ccc$ or simply by $\ccc$. If $\ddd \lra \ccc _S$ is an
isomorphism, the deformation is called \emph{strict}. If
$\ccc$ is flat over
$R$ (as an abelian category), the deformation is called \emph{flat}.
\end{definition}
\begin{remark}\label{ref-5.3-60}
In case of possible confusion we refer to the deformations of Definition
\ref{ref-5.1-58} as \emph{linear deformations} and to the deformations of
Definition \ref{ref-5.2-59} as \emph{abelian deformations}.
\end{remark}

We will now assume that $\theta:R\r S$ is a homomorphism between coherent
commutative rings such that $S$ is a finitely presented $R$-module. Then the
left and right adjoint to an abelian deformation
$\Dscr\r\Cscr$ exist and we will continue to denote them by $S\otimes_R-$ and
$\Hom_R(S,-)$.

Rephrazing the results of Sections \ref{ref-3-31},\ref{ref-4-40} in terms of
deformations we deduce that (flat) deformations are preserved under
some natural constructions.
\begin{proposition}
\label{ref-5.4-61}
Let $\BBB$ be an essentially small $S$-linear category.
An essentially small linear
       $R$-deformation $\AAA \lra \BBB$ induces an abelian deformation
       $\Mod(\mathfrak{b})\r \Mod(\mathfrak{a})$. Moreover flatness of these
       deformations is equivalent.
\end{proposition}
\begin{proof} We may assume that $\BBB=S \otimes_R\AAA$. Then the result
follows from Propositions \ref{ref-3.7-39} and \ref{ref-4.4-45}.
\end{proof}
\begin{proposition}\label{ref-5.5-62}
       Every abelian $R$-deformation $\ddd \lra \ccc$ between essentially
       small categories induces an abelian $R$-deformation $\Ind(\ddd) \lra
       \Ind(\ccc)$. Moreover flatness of these deformations is
       equivalent.
\end{proposition}
\begin{proof} We may assume $\Dscr=\Cscr_S$. Then the result follows
from Propositions \ref{ref-3.6-38} and \ref{ref-4.5-46}.
\end{proof}
\begin{proposition}\label{ref-5.6-63}
A deformation $\ddd \lra \ccc$ of locally coherent Grothendieck categories
induces a deformation $\Fp(\ddd) \lra \Fp(\ccc)$. Moreover flatness of
these deformations is equivalent.
\end{proposition}
\begin{proof} We may assume $\Dscr=\Cscr_S$. Then the result follows
from Propositions \ref{ref-3.6-38} and \ref{ref-4.6-49}(4).
\end{proof}

\section{Preservation of properties under nilpotent
deformations}\label{ref-6-64}

In this chapter we investigate some categorical/homological properties
that are preserved under flat, \emph{nilpotent} (see below) abelian
deformation.  More precisely, we will show that the following
properties of an abelian category $\ddd$ lift to a flat nilpotent deformation:
\begin{enumerate}
\item $\ddd$ is essentially small;
\item $\ddd$ has enough injectives;
\item $\ddd$ is a Grothendieck category;
\item $\ddd$ is a locally coherent Grothendieck category.
\end{enumerate}
lift to a flat nilpotent deformation.
\subsection{Nilpotent deformations}
As usual $\theta:R\r S$ is a morphism between coherent commutative rings
such that $S$ is finitely presented over $R$.
In order to lift properties from an $S$-linear (abelian) category to an
$R$-deformation, we obviously need some further assumptions on the ring map
$\theta$. First of all, we will assume that $\theta$ is
surjective. In this case, for an abelian $R$-linear $\ccc$, the forgetful
functor $\ccc_S \lra \ccc$ is fully faithful. Put $I=\ker \theta$. The
hypotheses imply that $I$ is finitely presented.

    Next, we will assume that the
kernel of $\theta$ is nilpotent.
\begin{definition}\label{ref-6.1-65}
The $R$-deformations
(linear or
abelian) are called \emph{nilpotent} provided the ideal $I$ is nilpotent, i.e.
$I^n = 0$ for some $n$. If $I^n = 0$
      the deformation is called
\emph{(nilpotent) of order $n$}.
\end{definition}
\begin{remark}
\label{ref-6.2-66}
If $\ccc$ a nilpotent deformation of $\ddd$ (linear or abelian),
then clearly there exist categories $\ddd = \ccc_0$, $\ccc_1$,
$\dots$, $\ccc_k = \ccc$ for certain $k$ such that $\Cscr_{i+1}$ is a
nilpotent deformation of order $2$ of $\Cscr_i$. If $\Cscr$ is flat
the we may assume that the $(\Cscr_i)_i$ are flat as well.
\end{remark}
\begin{proposition} (``Nakayama'')\label{ref-6.3-67}
Consider a nilpotent deformation $\Dscr \lra \ccc$ and $C \in \ccc$. If
either $S
\otimes_R C = 0$ or $\Hom_R(S,C) = 0$, then $C = 0$.
\begin{proof}
If $S\otimes_R C=0$ then $C=IC$ (\S \ref{ref-2.3-6}). Since $I$ is nilpotent
this implies
$C=0$. The case $\Hom_R(S,C) = 0$ is similar.
\end{proof}
\end{proposition}

\begin{proposition}\label{ref-6.4-68}
Consider a nilpotent $R$-deformation $\ddd \lra \ccc$. The functor
$\HomS: \ccc \lra \ddd$ reflects monomorphisms.
\begin{proof}
This readily follows from Proposition \ref{ref-6.3-67} using kernels.
\end{proof}
\end{proposition}

\subsection{Preservation of size}
In this section we temporarily drop the assumption that our
categories are automatically $\Uscr$-categories and that the base rings
$R,S$ are $\Uscr$-small.

We show that nilpotent deformations behave well with respect to size
matters.

\begin{lemma}\label{ref-6.5-69}
For $M \in \Mod(R)$, $$|M| \leq \sup\{|S \otimes_R M|,|\N|\}.$$
\begin{proof}
This follows by considering the $I$-adic filtration on $M$.
\end{proof}
\end{lemma}

For a category $\ccc$, the \emph{skeleton} $\Sk(\ccc)$ of $\ccc$ is the set
of all isomorphism classes of $\ccc$-objects.

\begin{lemma}\label{ref-6.6-70}
Consider an $S$-linear category $\BBB$ and let
$f:\AAA
\lra{\BBB}$ be a nilpotent $R$-deformation of $\BBB$.
Suppose $\alpha$ is an infinite cardinal such that
$|\Sk(\BBB)|\le \alpha$ and for all $B,B'\in \Ob(\BBB)$ one has
$|\BBB(B,B')|\le \alpha$. Then the analogous property holds for $\AAA$.
\begin{proof}
This follows from the following observations:
\begin{enumerate}
\item By Lemma \ref{ref-6.5-69}, $|\AAA(A,A')| \leq
     \sup\{|\BBB(f(A),f(A'))|,|\N|\}$ for $A, A' \in \AAA$;
\item By Proposition \ref{ref-A.1-141}, $|\Sk(\AAA)| = |\Sk(\BBB)|$.\hfill \qed
\end{enumerate}
\def\qed{}\end{proof}
\end{lemma}

\begin{proposition}
\label{ref-6.7-71}
Suppose $S$ is $\uuu$-small. Consider an $S$-linear category $\BBB$ and let
$f:\AAA
\lra{\BBB}$ be a nilpotent $R$-deformation of $\BBB$.
If $\BBB$ is a $\uuu$-category (resp. essentially $\uuu$-small),
the same holds for $\AAA$.
\begin{proof}
Immediate from (the proof of) Lemma \ref{ref-6.6-70}.
\end{proof}
\end{proposition}
We want to prove a similar result for abelian deformations. We prove
an analogue of Lemma \ref{ref-6.6-70}. Note that the result we prove
is more general than what we immediately need but its more
general form will be used afterwards.

\begin{lemma}
\label{ref-6.8-72}
Suppose $S$ is
$\uuu$-small. Consider a flat deformation
$\Dscr\r
\Cscr$ of abelian categories which is of order two. Let $\Dscr'$ be a Serre
subcategory of $\Dscr$ (i.e. a full subcategory closed under subquotients and
extensions) and let
$\Cscr'$ be the Serre subcategory
of $\Cscr$ generated by $\Dscr'$ (i.e. the full subcategory of
$\Cscr$ of objects $C$ with $S \otimes_R C$ (hence also $IC$) $\in \Dscr'$,
see also Proposition \ref{ref-7.5-111}).

Suppose $\alpha$ is an infinite cardinal such that
$|\Sk(\Dscr')|\le \alpha$ and for all $D,D'\in \Ob(\Dscr')$ one has
$|\Dscr(D,D')|\le \alpha$. Then the analogous property holds for $\Cscr'$.
\end{lemma}
\begin{proof}
Choose $C,C'$. Filtering $C,C'$ by the $I$-adic filtration
we immediately deduce
\begin{equation}
\label{ref-6.1-73}
|\ccc(C,C')| \leq |\ddd(S \otimes_R C, IC')|\cdot|\ddd(S \otimes_R C, S
\otimes_R C')|
\end{equation}
This implies $|\Cscr(C,C')|\le \alpha$

    Concerning
$|\Sk(\ccc')|$, we make the following observations:
\begin{enumerate}
\item An object in $\Cscr'$ is up to isomorphism determined by a triple
      $(A,B,e)$ where $A,B\in \Ob(\Dscr')$ and $e\in \Ext^1_\Cscr(A,B)$.
\item The spectral sequence \eqref{ref-4.1-51}
yields the bound:
      \[
      |\Ext^1_\Cscr(A,B)|\le |\Ext^1_\Dscr(A,B))|
      |\Hom_\Dscr(A,\Ext^1_R(S,B))|
      \]
\item An element of $\Ext^1_\Dscr(A,B)$ may be represented by a triple
      $(C,u,v)$ where $C\in \Ob(\Dscr')$ and $u\in \ddd(A,C)$, $v\in
      \ddd(C,B)$ are such that $(u,v)$ constitutes a short exact
      sequence.
      \end{enumerate}
The bound $|\Sk(\Cscr')|\le \alpha$ now easily follows.
      \end{proof}
\begin{proposition}\label{ref-6.9-74}
Suppose $S$ is $\uuu$-small. Consider a flat nilpotent deformation $\Dscr\r
\Cscr$ of abelian categories. If $\Dscr$ is a $\uuu$-category (resp. an
essentially $\uuu$-small category), then so is
$\Cscr$.
\end{proposition}
\begin{proof}
     We may assume that the deformation is of order two.  The case that
      $\Dscr$ is essentially $\uuu$-small follows directly from
     Lemma \ref{ref-6.8-72}. For the case that $\Dscr$ is a
     $\Uscr$-category we invoke \eqref{ref-6.1-73} to obtain that every
     $|\ccc(C,C')|$ is bounded by the cardinality of an element of $\uuu$.
\end{proof}
\subsection{Lifting of objects}\label{ref-6.3-75}

An important tool in the study of deformations is the lifting of
objects along the functors $(S \otimes_R -)$ and $\Hom_R(S,-)$. In
this section we state some results on lifting that will be used
afterwards. We start with the following definition.

\begin{definition}\label{ref-6.10-76}
     Consider a functor $H: \ccc \lra \ddd$ and an object $D$ of $\ddd$.
     A \emph{lift} of $D$ along $H$ is an object $C$ of $\ccc$ and an
     isomorphism $D \cong H(C)$.
       Lifts of $D$, together with the obvious
     morphisms, constitute a groupoid $L(H,D)$.
\end{definition}
\begin{remark}\label{ref-6.11-77}
If $H: \ccc \lra \ddd$ is right adjoint to a functor $F: \ddd \lra \ccc$, a
lift of $D$ along $H$ can be represented by a map $F(D) \lra C$. If $H$ is left
adjoint to $G: \ddd \lra \ccc$, a lift of $D$ along $H$ can be represented by a
map
$C \lra G(D)$.
\end{remark}
Consider a flat deformation $\ddd \lra \ccc$.
The restriction $\Hom_R^{cf}(S,-)$ of $\Hom_R(S,-)$ to coflat objects yields
a groupoid
$$L^{cf}(D) = L(\Hom_R^{cf}(S,-),D)$$
for every coflat $D$ in $\ddd$.
Analogously, the restriction $(S \otimes_R^{f} -)$ of $(S \otimes_R -)$ to
flat objects yields a groupoid
$$L^f(D) = L((S \otimes_R^{f} -),D)$$
for every flat $D$ in $\ddd$.

We state the following theorem without proof. The theorem is a special case
of an obstruction theory for derived lifting, which will be published
separately \cite{low2}.

\begin{theorem}\label{ref-6.12-78}
Consider a flat nilpotent $R$-deformation $\ddd \lra \ccc$ of order $2$ and a
coflat object
$D$ of $\ddd$. There is an obstruction $$o(D)
\in \Ext^2_{\ddd}(\Hom_S(I,D),D)$$ satisfying $$L^{cf}(D) \neq
\varnothing \iff o(D) = 0.$$ over
$$\Ext^1_{\ddd}(\Hom_S(I,D),D).$$
\end{theorem}

The following proposition shows that in a certain sense, the conditions of
Theorem \ref{ref-6.12-78} can itself be lifted under deformation.

\begin{proposition}\label{ref-6.13-79}
Consider a flat nilpotent deformation $\ddd \lra \ccc$ and coflat
$\ccc$-objects $\ovl{D}$ and $\ovl{E}$ with $\Hom_R(S,\ovl{D}) = D$ and
$\Hom_R(S,\ovl{E}) = E$. Suppose that for certain $i \geq 0$ we have that
$\Ext^i_{\ddd}(\Hom_S(X,D),E) = 0$ for all $X \in \modS$. Then
$\Ext^i_{\ccc}(\Hom_R(X,\ovl{D}),\ovl{E}) = 0$ for all $X \in \modR$.
\begin{proof}
We may consider $\ccc_S \lra \ccc$ of order $2$.
By Proposition \ref{ref-4.7-50},
$$\Ext^i_{\ccc}(\Hom_R(Y,\ovl{D}),\ovl{E}) = \Ext^i_{\ddd}(\Hom_S(Y,D),E) =
0$$ for all $Y \in \modS$.
For $X \in \modR$, it suffices to consider $0 \lra IX \lra X \lra S \otimes_R X
\lra 0$ and the corresponding $0 \lra \Hom_R(S \otimes_R X, \ovl{D}) \lra
\Hom_R(X,\ovl{D}) \lra \Hom_R(IX,\ovl{D}) \lra 0$. The result for $X$ then
follows from the result for $Y = IX$ and $Y = S \otimes_R X$ via the long
exact $\Ext$-sequence.
\end{proof}
\end{proposition}

\begin{corollary}\label{ref-6.14-80}
Consider a flat nilpotent deformation $\ddd \lra \ccc$ and a
coflat $\ddd$-object
$D$. If $\Ext^2_{\ddd}(\Hom_S(X,D),D) = 0$ for all $X \in \modS$, then
$\mathrm{Sk}(L^{cf}(D)) \neq \varnothing$. If
$\Ext^2_{\ddd}(\Hom_S(X,D),D) =
\Ext^1_{\ddd}(\Hom_S(X,D),D) = 0$, then $\mathrm{Sk}(L^{cf}(D))$
contains precisely one element.
\begin{proof}
This follows from Propositions \ref{ref-4.8-55}, \ref{ref-6.13-79} and
Theorem \ref{ref-6.12-78}.
\end{proof}
\end{corollary}

\begin{corollary}\label{ref-6.15-81}
Consider a flat nilpotent deformation $\ddd \lra \ccc$ and a $\ddd$-injective
$D$. Then $\mathrm{Sk}(L^{cf}(D))$ contains precisely one lift $D
\lra C$ and $C$ is a $\ccc$-injective.
\begin{proof}
We may consider $\ccc_S \lra \ccc$ of order 2.
The first part of the statement follows from Theorem \ref{ref-6.12-78}. To
prove that
$C$ is injective, first notice that by Proposition \ref{ref-4.7-50},
$\Ext^i_{\ccc}(A,C) = 0$
for all $A \in \ccc_S$. Finally, for $A \in \ccc$, it suffices to consider
the long exact $\Ext$ sequence associated to $0 \lra IA \lra A \lra S
\otimes_R A \lra 0$.
\end{proof}
\end{corollary}

\subsection{Deformations of categories with enough injectives} The results of
the previous section allow us to prove the following theorem:

\begin{theorem}\label{ref-6.16-82}
Consider a flat nilpotent deformation $\ddd \lra \ccc$. If $\ddd$ has
enough injectives then so does $\ccc$.
\begin{proof}
We may consider $\ccc_S \lra \ccc$ and an object $C$ in $\ccc$. Take a
$\ccc_S$-monomorphism $m: \Hom_R(S,C) \lra I$ to a $\ccc_S$-injective $I$ such
that there is a lift $\ovl{I}$ with $\Hom_R(S,\ovl{I}) = I$. By
Corollary \ref{ref-6.15-81}, $\ovl{I}$ is injective and so $m$ can be lifted
to a map
$\ovl{m}: C \lra \ovl{I}$ with $\Hom_R(S,\ovl{m}) = m$. By Proposition
\ref{ref-6.4-68}, $\ovl{m}$ is a $\ccc$-monomorphism.
\end{proof}
\end{theorem}

This result is false for non-flat nilpotent deformations, as the following
example shows.
\begin{example}
\label{ref-6.17-83}
       Let $k$ be a field and let $S=k$, $R=k[\epsilon]$, $\epsilon^2=0$.
       Let $\Dscr=\Mod(k)$. Let $\kappa$ be an infinite small cardinal and
       consider the following category:
\[
\Cscr=\{M\in \Mod(R)\mid \dim_k IM\le \kappa\}
\]
It is easy to see that $\Cscr$ is full abelian subcategory of $\Mod(R)$.
Clearly $\Dscr\r \Cscr$ is a nilpotent deformation of order two.

A skeletal subcategory of $\Mod(R)$ is given by the objects
\[
V(\alpha,\beta)=k^{\oplus\alpha} \bigoplus k[\epsilon]^{\oplus\beta}
\]
where $\alpha$, $\beta$ are small cardinals \cite{Ta}.
$\Cscr$ consists of those objects
$V(\alpha,\beta)$ for which $\beta\le \kappa$.

The objects $V(\alpha,\beta)$ with $\alpha\neq 0$ are not injective in
$\Cscr$ since we may always replace one copy of $k$ by $k[\epsilon]$
which yields a non-split extension.

Consider an object of the form
$V(\alpha,\kappa)$ in $\Cscr$ with $\alpha>\kappa$ and assume that
there is an injective object $E$ in $\Cscr$ containing
$V(\alpha,\kappa)$. Then clearly $E$ must be of the form
$V(\alpha',\kappa)$ with $\alpha'+\kappa\ge\alpha+\kappa$. In particular
$\alpha'>\kappa$. But then $E$ cannot be injective by the discussion
in the previous paragraph. This contradiction shows that $\Cscr$ does
not have enough injectives.
\end{example}

\subsection{Deformations of Grothendieck categories}\label{ref-6.5-84}

The aim of this section is to prove that the property of being a Grothendieck
category is preserved under flat nilpotent deformation.

For a small pre-additive category $\mathfrak{g}$ we write
$\Pre(\mathfrak{g})=
\Mod(\mathfrak{g}^{\op})$ for the category of (additive) presheaves
      on $\mathfrak{g}$. We recall the following version of the Gabriel-Popescu
theorem.
\begin{theorem}\label{ref-6.18-85}\cite{GP,Popescu,lowen1}
     Let
$\ddd$ be a Grothendieck category and let
$\GGGG$ be a small
     generating subcategory of $\ddd$. There is a localization (i.e. a fully
faithful functor with an exact left adjoint)
$$\ddd \lra\Pre(\GGGG): D \longmapsto \ddd(-,D)$$

Conversely if $\frak{g}$ is a small pre-additive category and $\Dscr\lra
\Pre(\GGGG)$ a localization, then $\ddd$ is a Grothendieck category.
Equivalence classes of such localizations are in one-one correspondence with
so-called Gabriel (or additive) topologies on
$\frak{g}$.
\end{theorem}

Let $\ddd$ be a Grothendieck category with a distinguished generator
$G$.  We define a \emph{cardinality} $|\cdot|$ on $\Ob(\ddd)$ by
$$|C| = |\ddd(G,C)|,$$ where
$|\ddd(G,C)|$ denotes the set theoretic cardinality of the set $\ddd(G,C)$.
For a cardinal $\kappa$, $\ddd_{\kappa}$ is the full subcategory of
$\ddd$ containing all $\ddd$-objects
       $D$ with $|D| \leq \kappa$.

\begin{example}
If we take $\ddd$ to be a module category $\Mod(A)$ over a ring $A$ and we
take $G = A$, the newly defined $|M|$ coincides with the set theoretic
cardinality of $M$.
\end{example}

The following proposition lists the essential properties of the
categories $\Dscr_\kappa$.
\begin{proposition}\label{ref-6.20-86}
     If $\kappa$ is a small cardinal, then $\ddd_{\kappa}$ is essentially
     small.  If $\kappa$ is an infinite cardinal of the form $2^\beta$,
     with $\beta\ge |G|$ then $\Dscr_\kappa$ is a Serre subcategory of
     $\Dscr$.
\end{proposition}
\begin{proof} \MMM We prove first that $\ddd_{\kappa}$ is essentially
     small.  By Theorem \ref{ref-6.18-85} $\Dscr_\kappa$ is equivalent to a full
subcategory of
     $\Mod(\Dscr(G,G)^{\op})_\kappa$. The latter obviously has small skeleton
     since the set of right $\ddd(G,G)$-module structures on every
     $\kappa'\le\kappa$ is small.

     The fact that $\Dscr_\kappa$ is a Serre subcategory of $\Dscr$ is a
     consequence of Lemmas \ref{ref-6.21-87},\ref{ref-6.22-88} below which show
     that the ordinary properties of cardinalities of modules more or
     less hold for objects in Grothendieck categories. The only not
     entirely obvious fact is that $\Dscr_\kappa$ is closed under
     quotients under the stated hypotheses on $\kappa$. This follows from
     Lemma \ref{ref-6.22-88} and the following computation.
\[
(2\kappa)^{|G|}=2^{\beta\cdot|G|}=2^\beta=\kappa\qed
\]
\def\qed{}\end{proof}
\begin{lemma}\label{ref-6.21-87}
For an
exact sequence
$0 \lra A \lra B \lra C$ in $\ddd$,
\begin{enumerate}
\item $|A| \leq |B|$;
\item $|B| \leq |A|\cdot|C|$.
\end{enumerate}
\begin{proof}
This follows from the exactness of $0 \lra \ddd(G,A) \lra \ddd(G,B) \lra
\ddd(G,C)$.
\end{proof}
\end{lemma}
\begin{lemma}\label{ref-6.22-88}
For an epimorphism $f: A \lra B$,$$|B| \leq (2|A|)^{|G|}.$$
\begin{proof}
     For a map $g: G \lra B$, consider the pullback $P_g$ of $f$ and $g$
     and then an epimorphism $\bigoplus_{i\in I} G\r P_g$.  The composition
$(f_i)_{i\in
       I}:\bigoplus_{i\in I} G\r P_g\r G$ is an epimorphism.  Put
$J=\{f_i\mid i\in I\}\subset \Dscr(G,G)$.

    In this way we still have an epimorphism
$(j)_{j\in J}:\bigoplus_{j\in J} G\r G$
which fits in a
    commutative diagram
$$\xymatrix{A \ar[r]^f & B \\ {\oplus_{j \in J}G} \ar[u]^{(a_j)_{j}}
\ar[r]_{(j)_j} & G\ar[u]_g}$$

Conversely it is clear that $g$ is uniquely determined by the data
$J\subset \Dscr(G,G)$ and $(a_j)_{j\in J}\in \Dscr(G,A)^J$ which yields
the required bound.
\end{proof}
\end{lemma}

In the remainder of this section, we consider a flat nilpotent deformation
$\Dscr=\Cscr_S
\lra
\ccc$ of order
$2$ in which $\Dscr$ is a Grothendieck category with a fixed generator $G$
defining a cardinality $|\cdot|$ on $\Dscr$. We put $\kappa=2^{\beta}$ for an
infinite cardinal $\beta \ge |G|$.
    By Proposition \ref{ref-6.20-86}, $\Dscr_\kappa$ is an essentially
small Serre subcategory of $\Dscr$. We let $\Cscr_\kappa$ be the
Serre subcategory of $\Cscr$ generated by $\Dscr_\kappa$. By
Lemma \ref{ref-6.8-72} $\Cscr_\kappa$ is essentially small as well.
Our aim is to show that $\Cscr_\kappa$ generates $\Cscr$.

The following lemma gives us a general procedure for constructing
generators for $\Cscr$. It is more general than what we need
but it will be reapplied in a slightly different setting
afterwards (see Theorem \ref{ref-6.36-103}).
\begin{lemma}\label{ref-6.23-89}
Let $(G_i)_{i\in I}$ be a (not necessarily small) collection of generators of
$\Dscr$. Consider a collection of objects $(G_f)_f$, indexed over all non-zero
$\ccc$-maps $f$ with codomain in $\Dscr$, obtained in the following manner:
\begin{equation}
\label{ref-6.2-90}
\xymatrix@=11pt{{G_f} \ar[rr]^-{k} \ar[rrdd]_{\gamma''}&&{G''_f} \ar[rr]^-{h'}
     \ar[dd]_{\gamma'} && {G'} \ar[rr]^-{1} \ar[dd]_{\gamma} && {G'}
     \ar[rr]^-{g'} \ar[dd]_{\beta'} && {A} \ar[rr]^-{f}
     \ar[dd]_{\beta} && {B} \\
&&&(P1)&&&&(P2)
\\
     &&{G_j} \ar[rr]_-{h} &&{S \otimes_R G'} \ar[rr]_-{\beta''} && {G_i}
     \ar[rr]_-{g} &&{S \otimes_R A} \ar[rruu]_{f'} &&}
\end{equation}
$g$ is chosen to make $f' \circ g \neq 0$, $G'$ is the pullback of
$\beta$ and $g$, $h$ is chosen to make $f' \circ g \circ \beta'' \circ h \neq
0$, $G''_f$ is the pullback of $\gamma$ and $h$, and $k$ is any map such
that $\gamma' =
\gamma \circ k$ is an epimorphism.
Then $(G_f)_f$ is a (not necessarily small) collection of generators for
$\ccc$.

\end{lemma}
\begin{proof}
     The $G_f$ clearly satisfy the generator property with respect to the
     maps $f:A\lra B$ with codomain in $\Dscr$. We claim this is
     sufficient.  Let $g:A\r C$ be a general map in $\Cscr$. If the
     composition $A\mathbin{\mathop{\longrightarrow}\limits^{g}} C\r S\otimes_R
C$ is not-zero then we
     are done.  If the composition is zero then we have a factorization
     $A\r IC\mathbin{\mathop{\longrightarrow}\limits^{}} C$
and we are done as well.
\end{proof}

\begin{lemma}\label{ref-6.24-91}
     Suppose $I=\{\ast\}, G_{\ast} = G$ and let $(G_f)_f$ be constructed
     as in Lemma \ref{ref-6.23-89} with every $k$ taken to be $1_{G''_f}$.
Then for all $f$ we have  $G_f \in \ccc_{\kappa}$.
    In particular, $\ccc$ has a small set of generators.
\end{lemma}
\begin{proof}
Recall that the squares marked with (P*) in \eqref{ref-6.2-90} are pullbacks.
  From (P1) it follows that
we have to show that $IG'\in \Dscr_\kappa$. From (P2) we obtain an exact
sequence
\[
0\lra IA \lra G'\r  G\lra 0
\]
Tensoring this sequence with $S$ we obtain a commutative diagram with
exact rows and columns
\[
\xymatrix{
&&0&&\\
\Tor_1^R(S,G) \ar[r] & IA \ar[r] & S\otimes_R G'\ar[r]\ar[u] &
G \ar[r] & 0 \\
0 \ar [r] & IA\ar[r]\ar@{=}[u] & G'\ar[r]\ar[u] & G\ar[r]\ar@{=}[u] & 0\\
&& IG'\ar[u] &&\\
&&0\ar[u]&&
}
\]
Diagram chasing yields an epimorphism
$
\Tor_1^R(S,G)\r IG'
$
which finishes the proof.
\end{proof}

\begin{proposition}\label{ref-6.25-92}
     Let $\Escr$ be an abelian category with enough
     injectives. Put $\III=\Inj(\Escr)$.
     Consider the functor
$$\Psi:
\Escr\op
\lra
\Mod(\III): C
\longmapsto \Escr(C,-).$$
\begin{enumerate}
\item $\Psi$ is exact, limit preserving and fully faithful;
\item
$\Psi$ induces an equivalence $$\Psi': \Escr\op \lra \mmod(\III)$$
\end{enumerate}
\end{proposition}
\begin{proof}
     $\Psi$ is exact and limit preserving since all the functors
     $\Escr(-,I): \Escr\op \lra \Ab$ are exact and limit preserving.

To show that $\Psi$ is fully faithful take $C,D \in \Ob(\Escr)$.
Choose an injective copresentation
\[
0\lra D\lra I_0\lra I_1
\]
Then we obtain a projective presentation of $\Escr(D,-)$ in $\Mod(\III)$
\[
\III(I_1,-)\r \III(I_0,-)\r \Escr(D,-)\r 0
\]
and thus
\begin{align*}
\Mod(\III)(\Escr(D,-),&\Escr(C,-))\\
&=
\ker (\Mod(\III)(\III(I_0,-),\Escr(C,-))\r
\Mod(\III)(\III(I_1,-),\Escr(C,-)))\\
&=\ker (\Escr(C,I_0)\r \Escr(C,I_1))\\
&=\Escr(C,D)
\end{align*}
    Now consider $F \in \mmod(\III)$. Then
$F$ has a presentation
$$\oplus_{i=1}^n \III(I_i,-) \lra
\oplus_{j=1}^m \III(I_j,-) \lra F \lra 0,$$ hence it follows that $F$ is
isomorphic to $\Escr(C,-)$ with $C$ defined by $0 \lra C \lra
\oplus_{j=1}^mI_j \lra \oplus_{i=1}^nI_i$. This proves that
$\Psi'$ is essentially surjective.
\end{proof}

\begin{lemma}\label{ref-6.26-93}
$\ccc$ has arbitrary small coproducts and they are exact.
\begin{proof}
Besides our fixed universe $\Uscr$ we introduce a larger universe $\Vscr$ such
that $\ccc$ is $\Vscr$-small.

By Theorem \ref{ref-6.16-82}, $\ccc$ has enough injectives and putting
$\III = \Inj(\ccc)$, by Proposition \ref{ref-6.25-92} we have an exact,
limit preserving functor $\Psi: \ccc\op \lra \Vscr{-}\Mod(\III): C
\longmapsto \ccc(C,-)$ inducing an equivalence $\Psi': \ccc\op \lra
\Vscr{-}\mmod(\III)$. We are to prove that $\ccc\op$ has exact
$\Uscr$-small products. Since $\Vscr{-}\Mod(\III)$ has even exact
$\Vscr$-small   products, it suffices to show
that for a $\Uscr$-small set of $\ccc$-objects $(C_i)_i$ the product
$\prod_i\ccc(C_i,-)$ in $\Vscr{-}\Mod(\III)$ is finitely presented.  Consider
the short exact sequences in $\ccc$
$$0 \lra IC_i \lra C_i \lra S \otimes_R C_i \lra 0.$$
In $\Vscr{-}\Mod(\III)$ we can take the product of the exact sequences
$$0 \lra \ccc(S \otimes_R C_i,-) \lra \ccc(C_i,-) \lra \ccc(IC_i,-) \lra 0$$
to obtain an exact sequence
$$0 \lra \ccc(\coprod_iS \otimes_R C_i,-) \lra \prod_i\ccc(C_i,-) \lra
\ccc(\coprod_iIC_i,-)
\lra 0,$$
where we used that $\ccc_S$ has small $\Uscr$-coproducts and $\ccc_S
\lra \ccc$ preserves
them.
Then $\prod_i\ccc(C_i,-)$ is finitely presented in
$\Vscr{-}\Mod(\III)$ as an extension
of finitely presented objects.
\end{proof}
\end{lemma}

The following proposition is the dual of \cite[A.3.2]{Neemanboek}.

\begin{proposition}\label{ref-6.27-94}
Consider an abelian category $\Escr$ with exact coproducts and any small
category $\iii$. The colimit functor
$$\colim: \Fun(\iii,\Escr) \lra \Escr: F \lra \colim(F)$$
is trivially right exact and its left derived functors
$L_i\colim$ exist. Moreover if $\Escr'$ is another abelian category with exact
coproducts and $\phi: \Escr \lra \Escr'$ is an exact functor preserving
coproducts, the following diagram commutes:
$$\xymatrix{{\Fun(\iii,\Escr)} \ar[d]_{L_i\colim_{\Escr}}
\ar[r]^-{(\phi \circ \cdot)} & {\Fun(\iii,\Escr')}
\ar[d]^{L_i\colim_{\Escr'}} \\
\Escr \ar[r]_-{\phi} & \Escr'.}$$
\end{proposition}
\begin{lemma}\label{ref-6.28-95}
$\ccc$ has exact filtered colimits.
\begin{proof}
By Lemma \ref{ref-6.26-93} we know that filtered colimits exist in $\ccc$, so
we are to prove that for a filtered category $\iii$, the functor
$$\colim: \Fun(\iii,\ccc) \lra \ccc: F \lra \colim(F)$$ is exact. By
Lemma \ref{ref-6.26-93} and Proposition \ref{ref-6.27-94}, it suffices to prove
that
$L_1\colim = 0$. Consider $F: \iii \lra \ccc$. Using the exact sequences $0
\lra IF(i)
\lra F(i) \lra S \otimes_R F(i) \lra 0$ and $\phi: \ccc_S \lra \ccc$, we
obtain an exact sequence $$0 \lra \phi \circ F' \lra F \lra \phi \circ F''
\lra 0$$ for functors $F'$ and $F'': \iii \lra \ccc_S$.
We obtain an exact sequence
$$L_1\colim(\phi \circ F') \lra L_1\colim(F) \lra L_1\colim(\phi \circ F''),$$
but since filtered colimits are exact in $\ccc_S$, it follows from
Proposition \ref{ref-6.27-94} that both ends of the sequence and hence also the
middle term $ L_1\colim(F)$ are zero, which proves our assertion.
\end{proof}
\end{lemma}

\begin{theorem}\label{ref-6.29-96}
Consider a flat nilpotent deformation $\ddd \lra \ccc$.

If $\ddd$ is a
Grothendieck category, the same holds for $\ccc$.
\begin{proof}
It suffices to consider a flat deformation $\ccc_S
\lra \ccc$ of order $2$. The result then immediately follows from
Lemmas \ref{ref-6.24-91}, \ref{ref-6.26-93} and \ref{ref-6.28-95}.
\end{proof}
\end{theorem}

\subsection{Deformations of locally coherent Grothendieck categories}

In this section we show that local coherence of a Grothendieck category is
also preserved under flat nilpotent deformation. We begin with some
preliminary results.

\begin{definition}\label{ref-6.30-97}\cite{Krause1}
An object $C$ in an abelian category $\ccc$ is called \emph{fp-injective}
provided that a monomorphism $X \lra Y$ between finitely presented objects
induces an epimorphism $\ccc(Y,C) \lra \ccc(X,C)$.
\end{definition}

\begin{proposition}\label{ref-6.31-98}
A filtered colimit of fp-injectives is fp-injective.
\begin{proof}
Let $\colim_iC_i$ be a filtered colimit of fp-injectives and $0 \lra X \lra
Y \lra Z \lra 0$ an exact sequence of finitely presented objects. Then all
sequences
$0 \lra \ccc(Z,C_i) \lra \ccc(Y,C_i) \lra \ccc(X,C_i) \lra 0$ are
exact hence so is their filtered colimit $0 \lra \ccc(Z,\colim_iC_i) \lra
\ccc(Y,\colim_iC_i) \lra \ccc(X,\colim_iC_i) \lra 0$.
\end{proof}
\end{proposition}

\begin{proposition}\label{ref-6.32-99} \cite[A.1]{Krause1}
In a locally
coherent Grothendieck category $\ccc$, the following are equivalent for $C \in
\ccc$:
\begin{enumerate}
\item $C$ is fp-injective;
\item $\Ext^1(Z,C) = 0$ for every $Z$ in $\Fp(\ccc)$;
\item $\Ext^i(Z,C) = 0$ for all $i > 0$ for every $Z$ in $\Fp(\ccc)$.\hfill\qed
\end{enumerate}
\end{proposition}

\begin{lemma}\label{ref-6.33-100}
Let $C$ be finitely presented in an abelian category $\ccc$
with enough injectives.
If filtered colimits of injectives are $\ccc(C,-)$-acyclic, then
$\Ext^n_{\ccc}(C,-)$ preserves filtered colimits for $n \geq 0$.
\begin{proof}
Consider a filtered colimit $\colim_iA_i$. By \cite{Groth1}, we can make a
functorial choice of injective resolutions $I_i^{\cdot}$ of the objects $A_i$.
Then $\colim_iI_i^{\cdot}$ is a
$\ccc(C,-)$-acyclic resolution of $\colim_iA_i$, hence we may compute
\begin{align*}
\Ext^n_{\ccc}(C,\colim_iA_i) & = H^n(\ccc(C,\colim_iI_i^{\cdot})) \\
& = H^n(\colim_i\ccc(C,I_i^{\cdot})) \\
& = \colim_i\Ext^n_{\ccc}(C,A_i).\hfill \qed\end{align*}
\def\qed{}\end{proof}
\end{lemma}

\begin{proposition}\label{ref-6.34-101}
Let $C$ be finitely presented in a locally coherent
Grothendieck category $\ccc$. Then $\Ext^n_{\ccc}(C,-)$ preserves filtered
colimits for $n \geq 0$.
\begin{proof}
This follows from Proposition \ref{ref-6.32-99} and Lemma \ref{ref-6.33-100}.
\end{proof}
\end{proposition}

\begin{lemma}\label{ref-6.35-102}
Consider a flat nilpotent deformation of Grothendieck categories $\ccc_S \lra
\ccc$ with $\ccc_S$ locally coherent. For a finitely presented object $C$ of
$\ccc_S$, $\Ext^n_{\ccc}(C,-)$ preserves filtered colimits for $n \geq 0$.
\begin{proof}
By Lemma \ref{ref-6.33-100} it suffices to show that for a filtered colimit
$\colim_iI_i$ of injectives
$\Ext^n_{\ccc}(C,\colim_iI_i) = 0$. Now by Corollary
\ref{ref-2.8-11}, $\colim_iI_i$ is coflat. Consequently, we may use
Propositions \ref{ref-4.7-50} and \ref{ref-6.34-101} to compute that
\begin{align*}
\Ext^n_{\ccc}(C,\colim_iI_i) & = \Ext^n_{\ccc_S}(C,\Hom_R(S,\colim_iI_i)) \\
& = \Ext^n_{\ccc_S}(C,\colim_i\Hom_R(S,I_i)) \\
& = \colim_i\Ext^n_{\ccc_S}(C,\Hom_R(S,I_i)) \\
& = \colim_i\Ext^n_{\ccc}(C,I_i) = 0.\hfill \qed
\end{align*}
\def\qed{}\end{proof}
\end{lemma}

We are now able to prove the main theorem.

\begin{theorem}\label{ref-6.36-103}
Consider a flat nilpotent deformation $\ddd \lra \ccc$ of Grothendieck
categories. If $\ddd$ is locally coherent, the same holds for $\ccc$.
\begin{proof}
We may consider $\ccc_S \lra \ccc$ of order $2$. Let
$(G_i)_i$ be a small set of coherent generators of $\ccc_S$. We want
to carry out the
construction of
Lemma \ref{ref-6.23-89} for specific maps $k$. Suppose we have constructed
the objects $G''_f$. Since $\ccc_S$ is locally coherent, $IG'$ is a
small filtered colimit $IG' = \colim_jF_j$ of coherent objects $F_j$ of
$\ccc_S$. Every canonical map $F_i \lra IG'$ induces a map
$\Ext^1_{\ccc}(G_j,F_i) \lra \Ext^1_{\ccc}(G_j,IG')$ mapping an extension $0
\lra F_i \lra E \lra G_j \lra 0$ to the pullback
$$\xymatrix{F_i \ar[d] \ar[r]& IG' \ar[d] \\ E \ar[r]\ar[rd] & P \ar[d]
\\&G_j}$$ By Lemma \ref{ref-6.35-102}, $\Ext^1_{\ccc}(G_j,-)$ preserves filtered
colimits, hence
$$\Ext^1_{\ccc}(G_j,IG') = \colim_i\Ext^1_{\ccc}(G_j,F_i).$$
It then follows from the construction of this filtered colimit that the
extension
$$\xymatrix{0 \ar[r]&IG' \ar[r]_--{\delta'} &G''_f \ar[r]_--{\gamma'} &G_j
\ar[r]&0}$$
is in the image of one of the maps $\Ext^1_{\ccc}(G_j,F_i) \lra
\Ext^1_{\ccc}(G_j,IG')$. Hence we obtain an extension
$$\xymatrix{0 \ar[r]&F_i \ar[r]_--{\delta''} &G_f \ar[r]_--{\gamma''} &G_j
\ar[r]&0}$$ and a map $k: G_f \lra G''_f$ as in Lemma \ref{ref-6.23-89}.
Hence $(G_f)_f$ is a set of generators for $\ccc$.
Now (see \S\ref{ref-2.2-5}) $\Fp(\ccc_S)$ has a small skeleton and since
$\ccc$ has small Yoneda Ext groups, it follows that the image of $(G_f)_f$ is
small. Using Proposition \ref{ref-4.6-47}, it is readily seen that the
objects $G_i$ and $F_i$ are coherent in $\ccc$. It then follows that all
objects of
$(G_f)_f$ are finitely generated, hence $\ccc$ is locally finitely generated.
But then all objects of $(G_f)_f$ are finitely presented, hence $\ccc$ is
locally finitely presented. Finaly all objects of $(G_f)_f$ are coherent hence
it follows that
$\ccc$ is locally coherent, as we set out to prove.
\end{proof}
\end{theorem}

\section{Deformation and localization}\label{ref-7-104}
\subsection{Statement of the main results}
\label{ref-7.1-105}
A fully faithful functor $i: \Lscr \lra \ccc$ between Grothendieck
categories is called a \emph{localization} if it has an exact left
adjoint. If in addition $i$ is an embedding of a full subcategory
closed under isomorphisms then we call $i$ a strict localization. In this
section, we study the compatibility of localization with deformations.

Below we prove the following result.
\begin{theorem} \label{ref-7.1-106}
Let $\Cscr_S\longrightarrow \Cscr$ be a nilpotent deformation
of Grothendieck categories.
      Then
$$\LLL \longmapsto \LLL_S$$
induces a bijection between the sets of strict localizations
of $\Cscr$ and $\Cscr_S$.
\end{theorem}
Note that this theorem \emph{does not} say that any nilpotent deformation of
a localization of $\Cscr_S$ is itself a localization of $\Cscr$. This
is in fact false (even under the appropriate flatness hypotheses).
\begin{example} Consider $S=k$, a field. $R=k[\epsilon]$, $\epsilon^2=0$.
      $\Lambda=k[x,y]$. $\Cscr=\Mod(\Lambda)$, $\Lscr=\Mod(\Lambda_x)$.
By Theorem \ref{ref-8.16-133} below the deformations of $\Cscr$ and $\Lscr$ correspond
to the deformations of $\Lambda$ and $\Lambda_x$. These are respectively
given by the rings
\begin{gather*}
R[x,y]/(yx-xy-f\epsilon)\qquad \text{and}\\
R[x,x^{-1},y]/((yx-xy-g\epsilon)
\end{gather*}
for $f\in k[x,y]$, $g\in k[x,x^{-1},y]$. Clearly $\Lambda_x$ has many
more deformations than $\Lambda$.
\end{example}
The next theorem allows us to recognize those deformations of
localizations which are themselves localizations.
\begin{theorem}\label{ref-7.3-107}
Consider a commutative diagram
\begin{equation}
\label{ref-7.1-108}
\xymatrix{{\ccc} \ar@<+2pt>[r]^{a'} &
{\LLL}\ar@<+2pt>[l]^{i'}\\ {\ccc_S} \ar[u]^F
\ar@<+2pt>[r]^a & {\LLL_S}
\ar[u]_{F'} \ar@<+2pt>[l]^i}
\end{equation}
in which $F$ and $F'$ are flat nilpotent deformations of Grothendieck
categories and
$(a,i)$ and $(a',i')$ are pairs of adjoint functors.
Suppose $i'$ maps injectives onto coflats.
If $i$ is a localization then the same holds for $i'$.
\end{theorem}
\WW If $u:\mathfrak{u}\r \Cscr$ is an additive functor from a small
pre-additive category $\UUU$ to a Grothendieck category $\ccc$, then
we have
       a pair of adjoint functors $(a,i)$ where
      $i:\Cscr\r
\Pre(\mathfrak{u}):U\longmapsto \Cscr(u(-),U)$ and $a:\Pre(\mathfrak{u})\r
\Cscr$ is the unique colimit preserving functor sending $\mathfrak{u}(-,A)$
to $u(A)$. We say that $u$ induces a localization if $i$ is a localization.

\medskip

We will apply  Theorem \ref{ref-7.3-107} in the following setting:
assume that $R\r S$ has nilpotent kernel and consider a commutative
diagram
\begin{equation}
\label{ref-7.2-109}
\xymatrix{{\UUU} \ar[r]^-u \ar[d]_f & {\ccc} \ar[d]^{(S \otimes_R -)}\\
{\VVV} \ar[r]_-v & {\ddd}}
\end{equation}
where the categories and functors are as follows:
\begin{enumerate}
\item $\mathfrak{u}$, $\mathfrak{v}$ are respectively small $R$ and $S$-linear
      flat categories;
\item the functor $f$ is a nilpotent linear deformation;
\item $\Cscr$, $\Dscr$ are respectively $R$ and $S$-linear
      flat  Grothendieck categories;
\item $u,v$ are respectively $R$ and $S$-linear additive functors;
\item the functor $S\otimes_R-$ is the left adjoint to a deformation
       $\ddd\r \ccc$;
\item the images of $u,v$ consist of flat objects.
\end{enumerate}
\begin{proposition}
\label{ref-7.4-110}
If $v$ induces a localization $\ddd\r
\Pre(\mathfrak{v})$, then $u$ induces a
localization $\ccc\r \Pre(\mathfrak{u})$.
\end{proposition}

\subsection{Proofs}
Let $\Cscr$ be an abelian category. Recall that a full subcategory
$\Sscr\subset\Cscr$ closed under subquotients and extensions in called a
\emph{Serre subcategory} of $\ccc$. If moreover $\ccc$ is Grothendieck and
$\sss$ is closed under coproducts, then
$\sss$ is called a
\emph{localizing subcategory} of $\Cscr$. If $\Sscr$ is localizing
then the quotient category $\Cscr/\Sscr$ is again a Grothendieck category
and the composition $\Sscr^{\perp}\r \Cscr\r\Cscr/\Sscr$, where
$\Sscr^\perp$ is the full subcategory of $\Sscr$ with objects
\[
\Ob(\sss^{\perp}) = \{ C \in \ccc\mid \forall
S
\in
\sss:\Ext^i_{\ccc}(S,C) = 0,\,\, i=0,1\},
\]
      is an equivalence of categories.

Assume that  $i: \Lscr \lra \ccc$ is a localization of Grothendieck
categories
with
exact left adjoint $a$.  Then $\ker(a)$ is a localizing
subcategory and $a$ induces an equivalence $\ccc/\ker(a) \lra \Lscr$.

If $i$ is a strict localization then $\Lscr = \ker(a)^{\perp}$. We
summarize this discussion in the following proposition:
\begin{proposition}\label{ref-7.5-111}
       Consider a Grothendieck category $\ccc$. Let $\SSS$ denote the set
       of localizing subcategories of $\ccc$ and let $\LLLL$ denote the
       set of strict localizations of $\ccc$. Then
$$\SSS \lra \LLLL: \sss \longmapsto \sss^{\perp}$$
and $$\LLLL \lra \SSS: (\LLL \xrightarrow{i} \ccc) \longmapsto
\ker(a),$$ where $a$ is a left adjoint of $i$, are inverse bijections.
\hfill \qed
\end{proposition}
For a set of objects or a subcategory $\Sscr\subset \Cscr$ we will
write $\langle\Sscr\rangle_{\Cscr}$ for the smallest Serre
subcategory of $\Cscr$ containing $\Sscr$.

Now assume that $R,S,I$ have their usual meaning and assume $I^n=0$.  Let
$\Cscr$ be an $R$-linear abelian category.
\begin{proposition}
\label{ref-7.6-112}
Let
$\SSS'$ denote the set of Serre subcategories of
$\ccc$ and $\SSS$ the set of Serre subcategories of $\ccc_S$. Then
\[
\SSS' \lra \SSS: \sss' \longmapsto \sss' \cap \ccc_S
\] and
\[
\SSS \lra \SSS': \sss \longmapsto \langle\sss\rangle_{\ccc}
\] are inverse bijections. If $\ccc$ is a Grothendieck category, they restrict
to inverse bijections between the respective sets of localizing subcategories.
\end{proposition}
\begin{proof} Let $\Sscr,\Sscr'$ be objects in $\mathfrak{s}$ and
$\mathfrak{s'}$
respectively.
Any object $C$ in $\Cscr$ has a finite filtration
\[
0 = I^nC \subseteq I^{n-1}C \subseteq
\dots
\subseteq I^2C
\subseteq IC \subseteq I^0C = C
\]
Write $\gr C\in \Cscr_S$ for its associated graded object. The formation
of $\gr C$ is compatible with coproducts.

For $\Sscr\in \mathfrak{s}$ let $\ovl{\Sscr}$ be the full subcategory of
$\Cscr$ whose objects are given  by
\[
\Ob(\ovl{\Sscr})=\{C\in \Cscr\mid \gr C\in \Sscr\}
\]
It is easy to see that $\ovl{\Sscr}$ is a Serre subcategory of
$\Cscr$ which is localizing if $\sss$ is.

Clearly $\Sscr\subset\ovl{\Sscr}\subset \langle \Sscr\rangle_\Cscr$ and since
$\ovl{\Sscr}$ is Serre, we deduce $\langle
\Sscr\rangle_\Cscr=\ovl{\Sscr}$.   This immediately implies
\[
\langle \Sscr\rangle_\Cscr\cap \Cscr_S=\Sscr
\]
Clearly if $\sss'\in
\mathfrak{s'}$ then
$\sss' \cap \ccc_S \in \SSS$, and since colimits in $\Cscr$ and $\Cscr_S$ are
computed in the same way (Proposition \ref{ref-4.2-42}), $\sss' \cap \ccc_S$ is
localizing if $\sss'$ is. We certainly have
\[
\langle \Sscr'\cap\Cscr_S\rangle_\Cscr\subset \Sscr'
\]
However if $C\in \Sscr'$ then $\gr C\in \Sscr'\cap\Cscr_S$ and hence
$C\in \langle \Sscr'\cap\Cscr_S\rangle_\Cscr$ by the earlier discussion.
Thus in fact $\langle \Sscr'\cap\Cscr_S\rangle_\Cscr= \Sscr'$ and
we are done.
\end{proof}
\begin{proof}[Proof of Theorem \ref{ref-7.1-106}]

       Let $i:\Lscr\r \Cscr$ be a strict localization and let $a:\Cscr\r
       \Lscr$ be the exact left adjoint to $i$. Consider the pair of
       functors $(i_S,a_S)$ between $\Lscr_S$ and $\Cscr_S$. It is clear
       that $i_S$ is still fully faithful and that $a_S$ is an exact left
       adjoint to $i_S$. Thus $i_S:\Lscr_S\r\Cscr_S$ is a (strict) localization.

Let $\LLLL'$ denote the set of strict localizations of $\ccc$
and $\LLLL$ the set of strict localizations of $\ccc_S$. Propositions
\ref{ref-7.5-111} and \ref{ref-7.6-112} furnish us with bijections
\[
\mathfrak{l}' \lra \mathfrak{s}'\lra \mathfrak{s} \lra \mathfrak{l}
\]
whose composition sends $\Lscr$ to $(\ker(a)\cap
\Cscr_S)^\perp=\ker(a_S)^\perp=\Lscr_S$, finishing the proof.
\end{proof}
The proof of Theorem \ref{ref-7.3-107} is based on the following observation:
\begin{proposition}\label{ref-7.7-113}
Consider a functor $i: \LLL \lra \ccc$ between Grothendieck categories
and suppose $i$ has a left adjoint $a$. The following are
equivalent:
\begin{enumerate}
\item $i$ is a localization;
\item $i$ preserves injectives and $i: \Inj(\LLL) \lra \Inj(\ccc)$ is fully
faithful.
\end{enumerate}
\end{proposition}
\begin{proof}
Obviously (1) implies (2). Suppose (2) holds. It is easily seen that $i$ is
fully faithful. The exactness of $a$ follows since exactness of a sequence can
be tested by considering $\Hom$'s into all injectives.
\end{proof}

\begin{proof}[Proof of Theorem \ref{ref-7.3-107}]
We use Proposition \ref{ref-7.7-113}. First, let us show that $i'$ preserves
injectives. Consider
$E
\in
\Inj(\LLL)$. To prove that $i'(E)$ is injective in $\ccc$, it suffices that
for all $C \in
\ccc_S$ we have $\Ext^1_{\ccc}(C,i'(E)) = \Ext^1_{\ccc_S}(C,\Hom_R(S,i'(E)))
=0$ where the first equality follows from Proposition \ref{ref-4.7-50}.

Looking at right adjoints in  \eqref{ref-7.1-108} we obtain a commutative
diagram
$$\xymatrix{
{\Lscr} \ar[r]^{i'} \ar[d]_{\Hom_R(S,-)} & {\Cscr} \ar[d]^{\Hom_R(S,-)}\\
{\Lscr_S} \ar[r]_i & {\Cscr_S.}}$$

The desired result follows from the fact that $\Hom_R(S,i'(E)) =
i(\Hom_R(S,E))$ is injective in $\ccc_S$. Next, we are to prove that
$\eta: \LLL(I,I') \lra \ccc(i'(I),i'(I'))$ is an isomorphism for
injectives $I$ and $I'$. Since $\LLL(I,I')$ and $\ccc(i'(I),i'(I'))$
are flat in $\Mod(R)$, using the 5-lemma it suffices that $S \otimes_R
\eta$ is an isomorphism. But $S \otimes_R \eta$ is isomorphic to
$\LLL_S(\Hom_R(S,I), \Hom_R(S,I')) \lra \ccc_S(i(\Hom_R(S,I)),
i(\Hom_R(S,I')))$, which proves our assertion.
\end{proof}
\begin{proof}[Proof of Proposition \ref{ref-7.4-110}]
      We appy Theorem \ref{ref-7.3-107}.
We have a diagram
\begin{equation}
\label{ref-7.3-114}
\xymatrix{%
{\Pre(\mathfrak{u})} \ar@<+2pt>[r]^{a'} &
{\ccc}\ar@<+2pt>[l]^{i'}\\
{\Pre(\mathfrak{v})} \ar[u]^F
\ar@<+2pt>[r]^a & {\ddd}
\ar[u]^{F'} \ar@<+2pt>[l]^i
}
\end{equation}
where $(a,i)$, $(a',i')$ are the pairs of adjoint functors associated
to $v$ and $u$; $F'$ is the deformation we
started with and $F$ is the abelian
deformation  associated to the linear deformation
$f^{\op}:\mathfrak{u}^{\op}\r\mathfrak{v}^{\op}$ (see Proposition \ref{ref-5.4-61}).

We claim that \eqref{ref-7.3-114} is commutative. We first consider the
``$i$-square''. Starting with $D\in \ddd$ we have
\begin{align*}
(F\circ i)(D)&=
\ddd((v\circ f)(-),D)\\
&=\ddd(S\otimes_R u(-),D)\\
&=\ccc(u(-),F'(D))\\
&=(i'\circ F')(D).
\end{align*}
Since $F,F'$ are fully faithful any $R$-linear functor $G$ between
$\Pre(\mathfrak{u})$ and $\ccc$ (in both directions) restricts to exactly
one functor between $\Pre(\mathfrak{v})$ and $\ddd$ (up to natural
isomorphism) necessarily given by $G_S$.  A pair of adjoint functors $(G,H)$
restricts to a pair of adjoint functors $(G_S,H_S)$.

Thus we obtain a pair of adjoint functors $(a'_S,i'_S)$ between
$\Pre(\mathfrak{u})$ and $\Cscr$ and furthermore $i=i'_S$. Since $a$ is
the left adjoint to $i$ we obtain $a=a'_S$ (up to natural isomorphism).
This implies the commutivity of the ``$a$-square'' in \eqref{ref-7.3-114}.

The remaining hypothesis of Theorem \ref{ref-7.3-107} we need to check is that
$i'$ sends injectives to coflats. If $E\in\ccc$ is injective then
it is coflat by Proposition \ref{ref-3.4-35}. Since the objects of
$\mathfrak{u}$ are mapped
to flat objects by $u$ this implies that $\ccc(u(-),E)$ takes on
coflat values by Proposition \ref{ref-2.9-12}(\ref{ref-7-21}). Hence by
Proposition \ref{ref-2.10-23}  $i'(E)$ is coflat.
\end{proof}

\section{Equivalent deformation problems}\label{ref-8-115}

\WW In this section we obtain several ``equivalences of deformation problems''.
We start by formalizing what we mean by this.

\subsection{Deformation pseudo functors}
In this section we need to be careful about
our choices of universe. Therefore we make them temporarily explicit
in our notations.

Let $\Uscr$ be a universe. We will denote by $\Uscr{-}\Rng^0$ the
category with as objects coherent commutative $\Uscr$-rings and as
morphisms surjective ring maps with a finitely generated, nilpotent
kernel. For a fixed coherent ring $S\in \Uscr$, we consider the category
$\Uscr{-}\Rng^0/S$.

Fix some other universe $\Wscr$. A deformation pseudo functor
is by definition a pseudo functor
\[
D:\Uscr{-}\Rng^0/S\r \Wscr{-}\Gd
\]
An \emph{equivalence} of deformation pseudo functors $D_1, D_2:
\Uscr{-}\Rng^0/S
\lra
\Wscr{-}\Gd$ is a pseudonatural transformation $\mu:
D_1 \lra D_2$ for which every $D_1(R) \lra D_2(R)$ is an equivalence of
categories. It is easy to see that this defines an equivalence relation
on deformation pseudo functors.

The dependence of our notations on the universes $\Uscr$ and $\Wscr$ is
a nuissance, but the deformation pseudo functors we will consider
below will be stable under enlarging $\uuu$ and $\www$ in a suitable sense.
The following proposition is a first
step in this direction.
\begin{proposition}
\label{ref-8.1-116}
Let $\Uscr\subset \uuu'$ be universes. Then
the obvious functor
\[
\Uscr{-}\Rng^0\!/S\r \uuu'{-}\Rng^0\!/S
\]
is an equivalence of categories. Consequently, two deformation pseudofunctors
$D_1, D_2: \Uscr'{-}\Rng^0/S
\lra \Wscr{-}\Gd$ are equivalent if and only if their restrictions to
$\uuu{-}\Rng^0/S$ are.
\end{proposition}
\begin{proof} We only need to show that it is essentially surjective, and
this readily follows from Proposition \ref{ref-6.5-69}.
\end{proof}

\subsection{Abelian and linear deformations}

\begin{definition}
\begin{enumerate} \item Consider an $S$-linear category $\BBB$.
An \emph{equivalence of deformations} from $f_1:\AAA_1 \lra
{\BBB}$ to $f_2:\AAA_2 \lra {\BBB}$ is an equivalence of
$R$-linear categories $\varphi:\AAA_2 \lra \AAA_1$ such that
$f_1 \circ \varphi$ is naturally isomorphic to $f_2$. If
$\varphi$ is an isomorphism with $f_1 \circ \varphi = f_2$, it is called an
\emph{isomorphism of deformations}.
\item Consider an abelian $S$-linear category $\ddd$.
An \emph{equivalence of deformations} from
$F_1: \ddd \lra \ccc_1$ to $F_2: \ddd \lra \ccc_2$ is an
equivalence of $R$-linear categories $\Phi: \ccc_1 \lra \ccc_2$ such that
$\Phi \circ F_1$ is naturally isomorphic to $F_2$. If $\Phi$ is an isomorphism
with $\Phi \circ F_1 = F_2$, it is called an
\emph{isomorphism of deformations}.
\end{enumerate}

\end{definition}
We consider the following groupoids for $R\in \Uscr{-}\Rng^0/S$ (where $S \in
\uuu$).
\begin{enumerate}
\item For a flat $S$-linear $\uuu$-category $\BBB$ and a universe $\vvv$ such
that
$\BBB$ is essentially $\vvv$-small and $\uuu \in \vvv$, we
       consider the groupoid $\vvv{-}\mathrm{def}_{\BBB}(R)$. The objects of
       $\vvv{-}\mathrm{def}_{\BBB}(R)$ are flat $R$-deformations of $\BBB$
which are
elements of $\vvv$.
       Its morphisms are equivalences of
       deformations modulo natural isomorphism of functors.
\item For a flat $S$-linear abelian $\uuu$-category $\ddd$ and a universe
$\vvv$ such that
$\ddd$ is essentially $\vvv$-small and $\uuu \in \vvv$, we
       consider the groupoid $\vvv{-}\mathrm{Def}_{\ddd}(R)$. The objects of
       $\vvv{-}\mathrm{Def}_{\ddd}(R)$ are the flat $R$-deformations of
       $\ddd$ which are elements of $\vvv$. Its morphisms are
       equivalences of deformations modulo natural isomorphism of functors.
\end{enumerate}

Clearly $\Ob(\vvv{-}\mathrm{def}_{\BBB}(R)) \subset \vvv$ and
$\Ob(\vvv{-}\mathrm{Def}_{\ddd}(R)) \subset \vvv$. Hence if we take a universe
$\www$ with $\vvv \in\Wscr$, this yields deformation pseudo functors
\[
\vvv{-}\deff_{\BBB}, \vvv{-}\Def_{\Dscr}:\Uscr{-}\Rng^0/\!S\r
\Wscr{-}\Gd
\]

The universe $\www$ is a purely technical device to make
sure $\vvv{-}\deff_{\BBB}$, $\vvv{-}\Def_{\Dscr}$ take their values in a
category. Obviously, the equivalence or non equivalence of two deformation
pseudofunctors $D_1, D_2:
\uuu{-}\Rng^0/\!S \lra \www{-}\Gd$ is preserved under changing from $\www$ to
$\www'$ for some $\www \subset \www'$.

The following proposition shows that the choice of the universe $\vvv$ is
harmless as well. Consider universes $\vvv \subset \vvv'$  with $\vvv' \in
\www$.

\begin{proposition}
\label{ref-8.3-117}
       Let $\mathfrak{b}$ be a flat $S$-linear $\Uscr$-category and let
$\Dscr$ be a flat $S$-linear abelian $\Uscr$-category.
Then the obvious pseudonatural tranformations
\begin{align*}
\vvv{-}\deff_{\mathfrak{b}}&\lra
\vvv'{-}\deff_{\mathfrak{b}}\\
\vvv{-}\Def_{\Dscr}&\lra \vvv'{-}\Def_{\Dscr}
\end{align*}
are equivalences of deformation pseudo functors.
\end{proposition}
\begin{proof} This follows from Lemmas \ref{ref-6.7-71},\ref{ref-6.9-74}.
\end{proof}

\subsection{Small skeletons}
  From now on, we will simply write $\ddef_{\BBB}$ and $\Def_{\ddd}$ to denote
the functors $\vvv{-}\ddef_{\BBB}$ and $\vvv{-}\Def_{\ddd}$ for some $\uuu$,
$\vvv$ and $\www$ as above.
Below we fix the universe $\Uscr$. We assume that all rings
are in $\Uscr$ and write $\Rng^0$ for $\Uscr{-}\Rng^0$. The notions of small
and essentially small are as usual with respect to the fixed universe $\Uscr$.
The same holds for the notion of a Grothendieck category.

We prove the following results.

\begin{theorem}\label{ref-8.4-118}
Assume that $\BBB$ is an essentially small flat $S$-linear category. Then
$\ddef_{\BBB}(R)$ has a small skeleton.
\end{theorem}
\begin{proof}
Since $\BBB$ is assentially small, there is an infinite small cardinal
$\kappa$ with $|\Sk(\BBB)| \le \kappa$ and $|\BBB(B,B')| \le \kappa$ for all
$B, B' \in \BBB$. By Lemma \ref{ref-6.6-70}, the same bound holds for
every $\AAA \in \ddef_{\BBB}(R)$. It is easily seen that
up to equivalence, the number of $R$-linear
categories $\AAA$ satisfying this bound is small. Furthermore, for a given
$\AAA$, up to natural isomorphism, the number of equivalences $S \otimes_R\AAA
\lra
\BBB$ is small (see Lemma \ref{ref-2.3-2}).
\end{proof}

\begin{theorem}
\label{ref-8.5-119}
     Assume that $\Dscr$ is a flat $S$-linear abelian category which
is either essentially small or
a Grothendieck category. Then $\Def_\Dscr(R)$ has a small skeleton.
\end{theorem}
\begin{proof}
We may assume that $I^2 = 0$. Assume first that $\ddd$ is essentially small.
Then there is an infinite small cardinal $\kappa$ with $|\Sk(\ddd)| \le \kappa$
and $|\ddd(D,D')| \le \kappa$ for all $D,D' \in \ddd$. By Lemma
\ref{ref-6.8-72}, the same bound holds for every $\ccc \in \Def_{\ddd}(R)$.
The proof is then finished like the proof of Theorem \ref{ref-8.4-118}.

Now assume that $\ddd$ is a Grothendieck category with a generator $G$.
Let $\kappa$ and $\ddd_{\kappa}$ be as in section
\ref{ref-6.5-84}. If $\ddd \lra \ccc$ is a deformation, the image $G'$ of $G$
in $\ccc_S$ defines a cardinality on $\ccc_S$ which allows us to define
$\ccc_{\kappa}$ as in section \ref{ref-6.5-84}. By Lemma \ref{ref-6.20-86},
we can take an infinite small cardinal $\lambda$ with
$|\Sk(\ddd_{\kappa})|
\le \lambda$ and
$|\ddd(D,D')| \le \lambda$ for all $D,D' \in \ddd_{\kappa}$. By Lemma
\ref{ref-6.8-72}, the same bound holds for $\ccc_{\kappa}$ for every $\ccc
\in
\Def_{\ddd}(R)$, and by Lemmas \ref{ref-6.23-89}, \ref{ref-6.24-91},
$\ccc_{\kappa}$ generates
$\ccc$. Consequently, using Theorem \ref{ref-6.18-85}, up to equivalence, a deformation
$\ddd
\lra
\ccc$ is determined by:
\begin{enumerate}
\item an $R$-linear category $\GGGG$ with $|\Ob(\GGGG)| \le
\lambda$ and $|\GGGG(G,G')| \le \lambda$ for all $G,G' \in \GGGG$;
\item an additive topology on $\GGGG$;
\item an object $G' \in \GGGG$ and a map $\ddd(G,G) \lra
\GGGG(G',G')$.
\end{enumerate}
Clearly, the number of such data
is small, which finishes the proof.
\end{proof}

\begin{remark}
\label{ref-8.6-120}
      Example \ref{ref-6.17-83} shows that Theorem \ref{ref-8.5-119} is
      trivially false without the flatness requirement on deformations, even for
      $\Dscr=\Mod(k)$, $k$ a field and $R=k[\epsilon]/(\epsilon^2)$.
      Indeed every small cardinal $\kappa$ in Example \ref{ref-6.17-83} yields a
      different deformation of $\Dscr$ and the set of small cardinals is
not itself small.
\end{remark}

\subsection{Elementary deformation equivalences}
We have the following elementary deformation equivalences.

\begin{proposition}\label{ref-8.7-121}
\begin{enumerate}
\item For every flat $S$-linear $\BBB$ there is an equivalence
$$\theta: \deff_{\BBB} \lra \deff_{\BBB^{\op}};$$
\item an equivalence $F: \BBB \lra \BBB'$ of flat $S$-linear
       categories induces an equivalence $$\delta: \deff_{\BBB} \lra
       \deff_{\BBB'};$$
\item for every flat abelian $\ddd$ there is an equivalence $$\theta:
\Def_{\ddd} \lra \Def_{\ddd^{\op}};$$
\item an equivalence $F: \ddd \lra \ddd'$ of flat abelian categories induces an
equivalence $$\delta: \Def_{\ddd'} \lra \Def_{\ddd'}.$$
\end{enumerate}
\end{proposition}
\begin{proof}
The only non-trivial point here is that by Proposition \ref{ref-3.3-34}, the
opposite of a flat abelian category is again flat. Note that the corresponding
statement for
$R$-linear categories is obviously true.
\end{proof}

\subsection{Ind-objects versus coherent objects}
\label{ref-8.5-122}
If $\ccc_{es}$ is essentially small then there is an equivalence
$\Cscr_{es}\cong\Fp(\Ind(\Cscr_{es}))$. Similarly for a locally coherent
Grothendieck category $\ccc$ there is an equivalence
$\Ind(\Fp(\ccc))\cong\ccc$ (see \S\ref{ref-2.2-4}).

Using Proposition \ref{ref-5.5-62} and Proposition \ref{ref-5.6-63}, we obtain pseudo
natural transformations
\begin{equation}
\label{ref-8.1-123}
\Def_{\Dscr_{es}}\r \Def_{\Ind(\Dscr_{es})}
\end{equation}
and
\begin{equation}
\label{ref-8.2-124}
\Def_{\Dscr}\r \Def_{\Fp(\Dscr)}
\end{equation}
for an essentially small $S$-linear category $\ddd_{es}$ and an $S$-linear
locally coherent Grothendieck category $\ddd$.
By Proposition \ref{ref-6.9-74}, any $\ccc \in \Def_{\ddd_{es}}$ is
essentially small, and by Theorem \ref{ref-6.36-103}, any $\ccc \in
\Def_{\ddd}$ is a locally coherent Grothendieck category. This allows one to
prove the following

\begin{theorem}\label{ref-8.8-125}
(\ref{ref-8.1-123}) and (\ref{ref-8.2-124}) are equivalences of deformation pseudo
functors. \hfill \qed
\end{theorem}

\subsection{Lifting deformations of localizations}
In this section we prove a general result  (Theorem \ref{ref-8.14-131} below)
which will allow us to construct
deformation equivalences in several settings afterwards.
Unfortunately the formulation of Theorem \ref{ref-8.14-131} is rather
technical so the interested reader may consider reading the subsequent
sections first.

\medskip

\let\ddef\deff Our purpose is to describe a class of functors $ v: \VVV
\lra \ddd$, with $\mathfrak{v}$ flat pre-additive and $\ddd$ flat abelian, for
which
$\ddef_{\VVV}$ and $\Def_{\ddd}$ are equivalent.

Consider a diagram
$$\xymatrix{{\UUU} \ar[d]_f & {\ccc} \ar[d]^{(S \otimes_R -)}\\
       {\VVV} \ar[r]_-u & {\ddd}}$$
in which $v$ is an $S$-linear functor
from a small $S$-linear category to an $S$-linear Grothendieck
category, $(S \otimes_R -)$ is left adjoint to an abelian
$R$-deformation $F: \ddd \lra \ccc$ and $f$ is a linear $R$-deformation.
In Section \ref{ref-7.1-105}, we have concentrated on lifting
properties of $v$, a lift of $v$ being given. Now we will concentrate
on the problem of lifting $v$. We have the following interpretation:

\begin{proposition}\label{ref-8.9-126}
The canonical functor
$$\xymatrix{{\Add(S)(\VVV, \ddd)} \ar[r]_-{\cong} & {\Add(R)(\UUU,\ddd)}
\ar[r]_-{(F
\circ \cdot)} & {\Add(R)(\UUU, \ccc)}}$$
is an abelian $R$-deformation, its left adjoint being given by
$$\xymatrix{{\Add(R)(\UUU, \ccc)} \ar[r]_-{(S \otimes_R) \circ \cdot)} &
{\Add(R)(\UUU,\ddd)} \ar[r]_-{\cong} & {\Add(S)(\VVV, \ddd).}}$$
Hence an $R$-linear $u: \UUU \lra \ccc$ making the diagram commute is a lift of
$v$ along this left adjoint.
\hfill \qed
\end{proposition}

\begin{proposition}\label{ref-8.10-127}
Consider an $S$-linear functor $v: \VVV \lra \ddd$ from a small $S$-linear
category to a
$S$-linear Grothendieck category such that the objects $v(V)$ are flat.
Suppose $\Ob(\VVV)$ is endowed with a reflexive, transitive relation
$\rrr$ such
that
\begin{enumerate}
\item $(U,V)$ not in
$\rrr$ implies $\VVV(U,V) = 0$;
\item $(U,V)$ in $\rrr$ implies that $v_{(U,V)}:
\VVV(U,V) \lra \ddd(v(U),v(V))$ is an isomorphism and that
$\Ext^i_{\ddd}(v(U), M \otimes_S v(V)) = 0$ for $M \in \mmod(S)$ and $i = 1,2$.
\end{enumerate}
Consider a flat nilpotent $R$-deformation $F: \ddd \lra \ccc$. The following
hold:
\begin{enumerate}
\item
There is a strict flat left
$R$-deformation $f: \UUU \lra \VVV$ and an $R$-linear functor $u: \UUU \lra
\ccc$ for which all objects $u(U)$ are flat such that the diagram
$$\xymatrix{{\UUU} \ar[r]^-u \ar[d]_f & {\ccc} \ar[d]^{(S \otimes_R -)}\\
{\VVV} \ar[r]_-v & {\ddd}}$$
commutes up to natural isomorphism.
\item If $f': \UUU' \lra \VVV$ and $u': \UUU' \lra \ccc$ are like $f$ and $u$
(but with $f'$ not necessarily strict), there is an equivalence $w: \UUU' \lra
\UUU$ with $f\circ w = f'$ and $v\circ w \cong v'$.
\end{enumerate}
\begin{proof}
We will constuct $u$ as stated in case $I^2 = 0$. It will be clear from
the proof and from Proposition \ref{ref-6.13-79} that the $R$-linear functor $u:
\UUU \lra \ccc$ satisfies the same properties as $v$, which then finishes the
proof in the general case.

For (1), consider a diagram
$$\xymatrix{ & {\ccc} \ar[d]^{(S \otimes_R -)}\\
{\VVV} \ar[r]_-v & {\ddd}}$$
as stated.
By Theorem \ref{ref-6.12-78}, $v(V)$ has a flat lift $D_V \lra v(V)$ along the
functor $(S
\otimes_R -)$. Define $\UUU$ by $\Ob(\UUU) = \Ob(\VVV)$ and $\UUU(U,V) =
\ccc(D_U,D_V)$ if $(U,V)$ is in $\rrr$ and $\UUU(U,V) = 0$ else. Since $\rrr$
is transitive, there is an obvious composition on $\UUU$ making $\UUU$ into
an $R$-linear category with an $R$-linear functor $u: \UUU \lra \ccc$.
There is also an obvious functor $f: \UUU \lra \VVV$. For $(U,V)$ not in
$\rrr$, $f_{(U,V)}: 0 = \UUU(U,V) \lra \VVV(U,V) = 0$ is trivially a
deformation of modules. For $(U,V)$ in $\rrr$, $f_{(U,V)}$ is isomorphic to
$(S \otimes_R -)_{(U,V)}: \ccc(D_U, D_V) \lra \ddd(S \otimes_R D_U, S
\otimes_R D_V)$. It follows from Proposition
\ref{ref-2.9-12}(\ref{ref-8-22}) that
$\VVV(U,V)$ is flat and $f_{(U,V)}$ is a deformation of modules. Hence $f$
is a flat left $R$-deformation, and the diagram commutes up to natural
isomorphism.

For (2), consider $f'$ and $u'$ as stated.
Since every
$u'(W)$ is flat and we have natural isomorphisms $S \otimes_R u'(W) \lra
v(f'(W))$, Theorem
\ref{ref-6.12-78} provides
$\ddd$-isomorphisms
$u'(W)
\lra D_{f'(W)}$ that allow the definition of a functor $w: \UUU' \lra \UUU$.
$w$ is essentially surjective since $f'$ is a deformation.
For
$(f'(U),f'(V))$ not in
$\rrr$, since
$f'_{(U,V)}:
\UUU'(U,V) \lra
\VVV(f'(U),f'(V)) = 0$ is a nilpotent deformation of modules, it follows that
$\UUU'(U,V) = 0$. Next, suppose $(f'(U),f'(V))$ is in $\rrr$.
We have a commutative diagram
$$\xymatrix{{\UUU'(U,V)} \ar[r]^-{u'_{(U,V)}} \ar[d]_{f'_{(U,V)}} &
{\ccc(u'(U),u'(V))}
\ar[d]\\ {\VVV(f'(U),f'(V))} \ar[r]_-{v_{(f'(U),f'(V))}} &
{\ddd(vf'(U),vf'(V))}}$$
in which the right arow is isomorphic to $(S
\otimes_R -)_{(u'(U),u'(V))}$ and hence both vertical arrows are flat
deformations of modules and
$v_{(f'(U),f'(V))}$ is an isomorphism. Since $I \otimes_R u'_{(U,V)} \cong I
\otimes_S v_{(f'(U),f'(V))}$, it follows from the 5-lemma that
$u'_{(U,V)}$ and hence also $w_{(U,V)}$ is an isomorphism. This proves that
$w$ is fully faithful.
\end{proof}
\end{proposition}

\begin{proposition}\label{ref-8.11-128}
For $v: \VVV \lra \ddd$ as in Proposition \ref{ref-8.10-127} and $R \in
\Rng^0/S$, there is a functor
$$\kappa_R: \Def_{\ddd}(R) \lra \ddef_{\VVV}(R).$$
\begin{proof}
For $F: \ddd \lra \ccc$ in $\Def_{\ddd}(R)$, we can define $\kappa_R(F)$
as $f: \UUU \lra \VVV$ of Proposition \ref{ref-8.10-127}(1). Proposition
\ref{ref-8.10-127}(2)
can be used to define $\kappa_R$ on morphisms.
\end{proof}
\end{proposition}

\begin{proposition}\label{ref-8.12-129}
For a localization of $S$-linear abelian categories
$$\xymatrix{{\ddd} \ar@<+2pt>[r]^a & {\kkk} \ar@<+2pt>[l]^i}$$
and $R \in \Rng^0/S$, there is a functor
$$\lambda_R: \Def_{\ddd}(R) \lra \Def_{\kkk}(R).$$
\begin{proof}
Since every localization factors as an
equivalence followed by a strict localization, this is a consequence of
Theorem \ref{ref-7.1-106} and Lemma \ref{ref-8.13-130} below applied to the functor
$a$.
\end{proof}
\end{proposition}

\begin{lemma}\label{ref-8.13-130}
Consider an $R$-linear exact functor $F: \ccc \lra \ddd$ that is essentially
surjective. If $\ccc$ is flat over $R$, the same holds for $\ddd$.
\begin{proof}
For $D \in \ddd$ and $X \in \mmod(R)$,
take an isomorphism $d: D \lra F(C)$ for a $C$ in $\ccc$.
Take a $\ccc$-monomorphism $c: C \lra C'$ with
$\Ext^1_{\ccc,R}(X,c) = 0$. Then $F(c) \circ d$ is a
$\ddd$-monomorphism and by Proposition \ref{ref-2.7-9},
$\Ext^1_{\ddd,R}(X,F(c)) = F(\Ext^1_{\ccc,R}(X,c)) = 0$.
\end{proof}
\end{lemma}

The construction of an abelian deformation of $\Mod(\BBB)$ from an $R$-linear
deformation of an essentially small category $\BBB$ yields a
pseudonatural transformation
$$\mu: \ddef_{\BBB} \lra \Def_{\Mod(\BBB)}.$$

\begin{theorem}\label{ref-8.14-131}
Let $v: \VVV \lra \ddd$ be as in Proposition \ref{ref-8.10-127} and suppose $v$
induces a localization
$$\xymatrix{{\Pre(\VVV)} \ar@<+2pt>[r]^a & {\ccc.} \ar@<+2pt>[l]^i}$$
The functors
$$\xymatrix{{\ddef_{\VVV}(R)} \ar[r]_-{\theta_R} & {\ddef_{\VVV\op}(R)}
\ar[r]_-{\mu_R} & {\Def_{\Pre(\VVV)}(R)} \ar[r]_-{\lambda_R} &
{\Def_{\ddd}(R)}}$$  and
$$\xymatrix{{\Def_{\ddd}(R)} \ar[r]_-{\kappa_R} & {\ddef_{\VVV}(R)}}$$
are inverse equivalences,
where $\lambda_R$ is as in Corollary \ref{ref-8.12-129}, $\kappa_R$ as
in Corollary \ref{ref-8.11-128},
$\theta_R$ as in Proposition \ref{ref-8.7-121} and $\mu_R$ as above.
\begin{proof}
First, we indicate a natural isomorphism $\lambda_R \circ \mu_R \circ \theta_R
\circ \kappa_R
\cong 1_{\Def_{\ddd}(R)}$. For $F: \ddd \lra \ccc$ in $\Def_{\ddd}(R)$, its
image
$F'$ under $\mu_R \circ \theta_R \circ \kappa_R$ fits into a diagram
$$\xymatrix{{\Pre(\UUU)} & {\ccc} \ar[l]_{i'} \\
{\Pre(\VVV)} \ar[u]^{F'} & {\ddd} \ar[u]_F \ar[l]_{i}}$$
that commutes up to natural isomorphism. By Proposition
\ref{ref-7.4-110}, $i'$ is
a localization. The construction of $\lambda_R(F')$ yields a diagram
$$\xymatrix{{\Pre(\UUU)} & {\ccc'} \ar[l]_{i'_1} &{\ccc} \ar[l]_{i'_0}^{\sim}\\
{\Pre(\UUU)_S} \ar[u] & {\ccc'_S} \ar[u] \ar[l]_{{i'_1}_S} &{\ccc_S}
\ar[l]_{\sim} \ar[u]\\
{\Pre(\VVV)} \ar[u]^{\sim} && {\ddd} \ar[ul]^{\sim} \ar[u]_{\sim} \ar[ll]^i}$$
in which $i'_1$ and ${i'_1}_S$ are strict localizations and $i'_0$ constitutes
an equivalence of deformations between $F$ and $\lambda_R(F')$.

Next, we indicate a natural isomorphism $\kappa_R \circ \lambda_R \circ \mu_R
\circ \theta_R \cong 1_{\ddef_{\VVV}(R)}$. Consider $f': \UUU' \lra \VVV$ in
$\ddef_{\VVV}(R)$. Its image
$F: \ddd \lra \ccc$ under $\lambda_R \circ \mu_R \circ \theta_R$ corresponds to
a localization $i': \ccc \lra \Pre(\UUU')$. The diagram
$$\xymatrix{{\UUU'} \ar[r]^-{u'} \ar[d]_{f'} & {\ccc} \ar[d]^{(S
\otimes_R -)}\\
{\VVV} \ar[r]_-v & {\ddd}}$$
in which $u'$ is induced by a left adjoint of $i'$ and $(S \otimes_R -)$ is a
left adjoint of $F$ commutes up to natural isomorphism. The functor $w: \UUU'
\lra \UUU$ of Proposition \ref{ref-8.10-127} constitutes an equivalence of
deformations between $f'$ and
$\kappa_R(F)$.
\end{proof}
\end{theorem}
\subsection{Deformations of Grothendieck categories with acyclic generators}

Next we will formulate some specializations of Theorem \ref{ref-8.14-131}.

\begin{theorem}\label{ref-8.15-132}
       Consider a flat $S$-linear Grothendieck category $\ddd$ and a small full
       subcategory $\GGGG$ such that $\Ob(\GGGG)$ is a collection of flat
       generators of $\ddd$ and such that for $G, G'$ in $\GGGG$
$$\Ext_{\ddd}^i(G, M \otimes_S G') = 0$$
for $i = 1,2$ and $M \in \mmod(S)$.
For $R \in \Rng^0/S$, there is an equivalence
$$\ddef_{\GGGG}(R) \cong \Def_{\ddd}(R).$$
\begin{proof}
It is readily seen that the inclusion $\GGGG \lra \ddd$ satisfies the
conditions of Theorem \ref{ref-8.14-131}.
\end{proof}
\end{theorem}

\begin{theorem}\label{ref-8.16-133}
For an essentially small flat $S$-linear category $\BBB$,
$$\mu: \ddef_{\BBB} \lra \Def_{\Mod(\BBB)}$$ is an equivalence of deformation
pseudo functors.
\begin{proof}
Since $\BBB\op$ is a full subcategory of $\Mod(\BBB)$
satisfying the conditions of Theorem \ref{ref-8.15-132}, the result follows from
Theorem \ref{ref-8.14-131}.
\end{proof}
\end{theorem}

\subsection{Deformations of categories with enough injectives}
\label{ref-8.8-134}

We have seen in Proposition \ref{ref-6.25-92} that if $\ccc$ is an
abelian category with enough injectives, there is an equivalence $\ccc
\cong (\mmod(\Inj(\ccc)))^{\op}$. We  prove the following

\begin{theorem} \label{ref-8.17-135} If $\Dscr$ has enough injectives
then $\Def_\Dscr$ and $\deff_{\Inj(\Dscr)}$ are equivalent deformation
pseudo functors.
\end{theorem}
\begin{proof}
Since $\mmod(\Inj(\ccc)) \cong \vvv{-}\Fp(\vvv{-}\Mod(\Inj(\ccc)))$ where
$\Inj(\ccc)
\in
\vvv$, we obtain a pseudo natural transformation
\begin{equation}
\label{ref-8.3-136}
\ddef_{\Inj(\ccc)} \lra \Def_{(\mmod(\Inj(\ccc)))^{\op}},
\end{equation}
which is actually an equivalence by Theorems \ref{ref-8.16-133}, \ref{ref-8.8-125} and
Proposition \ref{ref-8.7-121}.
\end{proof}

In Appendix A, we clarify this equivalence a little further using
certain preservation properties of flat nilpotent linear deformations.

\subsection{Sheaves of modules over a ringed space}
We will now give an application of Theorem \ref{ref-8.14-131}  to sheaves of
modules over a ringed space.

Let $X$ be a topological space and let $\ooo_X$
be a sheaf of $S$-algebra's on $X$.
For an open $U \subset X$,  Let
$\ooo_U$ be the restriction of $\ooo_X$ to $U$ and
denote by $\Mod(\ooo_U)$
and $\PMod(\ooo_U)$ the ($S$-linear) categories of sheaves
and presheaves of $\ooo_U$-modules.
For a basis $\bbb$ of the topology, denote by $\ooo_{\bbb}$ the restriction
of $\ooo_X$ to $\bbb$ and by $\PMod(\ooo_{\bbb})$ the category of presheaves
of $\ooo_{\bbb}$-modules.

\begin{theorem}
       Let $X$ be a topological space and let $\ooo_X$ be a sheaf of flat
       $S$-algebra's on $X$ as above. Suppose the topology has a basis
       $\bbb$ such that for every $U \in \bbb$, the sheaf cohomology
$$H^i(U,M \otimes_S \ooo_U) = 0$$
for $i = 1,2$ and $M \in \mmod(S)$.
For $R \in \Rng^0/S$, there is an equivalence
$$\Def_{\PMod(\ooo_{\bbb})}(R) \cong \Def_{\Mod(\ooo_{X})}(R).$$

\begin{proof}
Let $a: \PMod(\ooo_X) \lra \Mod(\ooo_X)$ be the exact
sheafication functor left adjoint to inclusion.
For an open $U \subset X$, let $j_U: U \lra X$ be the inclusion map and
Let $P_U = j^p_{U,!}\ooo_U$ and $S_U = j_{U,!}\ooo_U = a(j^p_{U,!}\ooo_U)$ be
the extensions by zero of $\ooo_U$ in the categories of presheaves and sheaves.
Also, let $P_U^b$ denote the restriction of $P_U$ to $\bbb$.
For any basis $\bbb$, for $U \in \bbb$ and $F \in \PMod(\ooo_{\bbb})$, we have
$\Hom(P_U^b,F) = F(U)$.
For $U$ and $V$ in $\bbb$, we obtain
\begin{equation}
\label{ref-8.4-137}
\Hom(P_U^b,P_V^b) = \Hom(P_U,P_V)=
\left\{
\begin{array}{r@{\hskip 1cm}l}
\ooo_X(U)& \mbox{if } U\subset V\\
0& \mbox{otherwise}
\end{array}
\right.
\end{equation}
and
\begin{equation}
\label{ref-8.5-138}
          \Hom(S_U,S_V)= \{f\in \ooo_X(U)\mid f {\mbox{ is zero on a
neighbourhood of }} U\setminus V\}.
\end{equation}
Thus in particular if $U\subset V$,
\begin{equation}
\label{ref-8.6-139}
\Hom(P_U,P_V)=\Hom(S_U,S_V).
\end{equation}
Let $\UUU$ be the full subcategory of $\PMod(\ooo_X)$ spanned by the objects
$P_U$ for $U \in \bbb$.
By (\ref{ref-8.4-137}), $\UUU$ is isomorphic to the full subcategory of
$\PMod(\ooo_{\bbb})$ spanned by the objects $P_U^b$ for $U \in \bbb$. Since
these objects are a family of finitely generated projective generators of
$\PMod(\ooo_{\bbb})$, we deduce that
$$\PMod(\ooo_{\bbb}) \cong \Pre(\UUU)$$
hence it suffices to prove that
$$\ddef_{\UUU}(R) \cong \Def_{\Mod(\ooo_{X})}(R).$$
Let $$u: \UUU \lra \Mod(\ooo_{X})$$
be the restriction of $a$.
We will show that $u$ satisfies the
conditions of Proposition \ref{ref-8.14-131}.
It is easily seen that $u$ induces a localization (see for example
\cite{lowen1}), so it remains to verify the conditions stated in
Proposition \ref{ref-8.10-127}. If we say that $(P_U,P_V)$ is in $\rrr$
if $U \subset
V$, the result follows from (\ref{ref-8.4-137}), (\ref{ref-8.6-139}) and the
following computation for $U \subset V$
$$\begin{aligned}
\Ext^i_{\Mod(\ooo_X)}(S_U, M \otimes_S S_V)
&= \Ext^i_{\Mod(\ooo_U)}(\ooo_U, (M \otimes_S S_V)|_U))\\
&= \Ext^i_{\Mod(\ooo_U)}(\ooo_U, M \otimes_S \ooo_U)\\
&= H^i(U, M \otimes_S \ooo_U)\\
\end{aligned}$$
where we have used that $j_{U,!}$ is exact and hence its right adjoint
restriction functor $(-)_U$ preserves injectives.
\end{proof}
\end{theorem}

\appendix
\section{\strut}
\label{ref-A-140}

In this Appendix we indicate a direct proof of Corollary \ref{ref-8.17-135}. This
proof makes use of some preservation properties of flat nilpotent linear
deformations that may be of independent interest.

\subsection{Flat nilpotent linear deformations}

let $F:\AAA \lra
\BBB$ be a flat nilpotent $R$-deformation of an $S$-linear category $\BBB$.
In this section, we will lift some properties of $\BBB$ to $\AAA$. The
following is well known in the ring case:

\begin{proposition}\label{ref-A.1-141}
If $f: A_1 \lra A_2$ in $\AAA$ is such that $F(f)$ is an isomorphism in $\BBB$,
then $f$ is an isomorphism in $\AAA$.
\begin{proof}
It suffices to consider $\AAA \lra S \otimes_R\AAA$ nilpotent of
order $2$. Suppose $f: A_1 \lra A_2$ in $\AAA$ is an isomorphism
in $\BBB$. Take an
$\AAA$-map $g: A_2 \lra A_1$ such that
$g \circ f = 1_{A_1}$ and $f \circ g = 1_{A_2}$ in $S \otimes_R \AAA$.
        From the exact rows $0 \lra I\AAA(A,A') \lra
\AAA(A,A') \lra (S \otimes_R\AAA)(A,A') \lra 0$ and from $I^2 = 0$
we deduce that $$(g \circ f - 1_{A_1})^2 =0$$ and
$$(f \circ g - 1_{A_2})^2 = 0.$$
These equations can be rewritten as
$$(2g - g\circ
f \circ g) \circ f = 1_{A_1}$$ and
$$f \circ (2g - g\circ
f \circ g) = 1_{A_2},$$
which proves our assertion.
\end{proof}
\end{proposition}

\begin{proposition}\label{ref-A.2-142}
If $Z \in \AAA$ is such that $F(Z)$
is a zero-object in $\BBB$, then $Z$ is a zero-object
in $\AAA$.
\begin{proof}
It suffices to consider $\AAA \lra S \otimes_R\AAA$. Suppose $Z$ is a
zero-object in $S \otimes_R\AAA$. So $(S \otimes_R\AAA)(A,Z) = 0$ for all
$\AAA$-objects $A$. But then by Nakayama $\AAA(A,Z) = 0$ for
all $A$, meaning that $Z$ is a zero-object in $\AAA$.
\end{proof}
\end{proposition}

\begin{definition}\label{ref-A.3-143}
In a category $\ccc$, an \emph{idempotent} is a map $e$ with $e \circ e = e$.
An idempotent $e$ \emph{splits} if there exist maps $r$,$s$ with $e = s \circ
r$ and $r \circ s = 1$. A category in which all idempotents split is called
\emph{Karoubian}.
\end{definition}

\begin{remark}\label{ref-A.4-144}
The splitting of an idempotent $e$ is equivalent to the existence
of the equalizer of $e$ and $1$ and to the existence of the coequalizer of $e$
and
$1$. Thus an abelian category is Karoubian. If $\ccc$ is an abelian
category with enough injectives, $\Inj(\ccc)$ is Karoubian too since a
retract of an injective is injective.
\end{remark}

\begin{proposition}\label{ref-A.5-145}
If $\BBB$ is Karoubian, the same
holds for $\AAA$.
\begin{proof}
It suffices to consider $\AAA \lra S \otimes_R\AAA$ nilpotent of
order $2$. Let $e: A \lra A$ be an idempotent in $\AAA$
and take maps $s: B \lra A$, $r: A \lra B$ in $\AAA$
such that $s \circ r = e$ and $r \circ s = 1$ in $S \otimes_R\AAA$.
        From $s \circ r - e \in I\AAA(A,A)$ and $r \circ s
- 1 \in I\AAA(B,B)$ we obtain the following equations
in $\AAA$:
\begin{enumerate}
\item $(s \circ r - e)^2 = 0$;
\item $(r \circ s - 1)^2 = 0$;
\item $(r \circ s - 1) \circ r \circ (s \circ r - e) = 0$;
\item $(s \circ r - e) \circ s \circ (r \circ s - 1)
= 0$.
\end{enumerate}
It follows from (1) that $$(s + e \circ s - s \circ r
\circ s) \circ (r \circ e) = e$$ and from (2) that
$$2(r \circ s) - r \circ s \circ r \circ s = 1.$$
Combining (1) and (3) gives us
$$r \circ e \circ s \circ r = r \circ s \circ r,$$
and combining this with (2) and (4) gives
$$r \circ e \circ s = r \circ s.$$
We can now compute that $(r \circ e) \circ (s + e \circ s - s \circ r
\circ s) = 2(r \circ s) - r \circ s \circ r \circ s
= 1$. So we have shown that $(s + e \circ s - s \circ r
\circ s): B \lra A$ and $(r \circ e): A \lra B$ give
a splitting of $e$ in $\AAA$.
\end{proof}
\end{proposition}

\begin{proposition}\label{ref-A.6-146}
If $A$, $B$, $C$, $s_1: A \lra C$ and $s_2: B \lra C$ in $\AAA$ are such that
$(F(C),F(s_1),F(s_2))$ is a coproduct of $F(A)$ and
$F(B)$ in $\BBB$, then $(C,s_1,s_2)$ is a coproduct
of $A$ and $B$ in $\AAA$.
\begin{proof}
It suffices to consider $\AAA \lra S \otimes_R\AAA$ nilpotent of
order $2$. Take maps $p_1: C \lra A$ and $p_2: C \lra
B$ in $\AAA$ such that $(C,s_1,s_2,p_1,p_2)$ is a biproduct
of $A$ and $B$ in $S \otimes_R\AAA$. We obtain the following
equations in $\AAA$:
\begin{enumerate}
\item $(p_1 \circ s_1 - 1)^2 = 0$;
\item $(p_1 \circ s_2) \circ (p_2 \circ s_1) = 0$;
\item $(p_1 \circ s_1 - 1) \circ (p_1 \circ s_2) =
0$;
\item $(p_1 \circ s_2) \circ (p_2 \circ s_2 - 1) =
0$;
\item $(s_1 \circ p_1 + s_2 \circ p_2 - 1)^2 = 0$.
\end{enumerate}
Put $$p_1' = 2p_1 - p_1 \circ s_1 \circ p_1
- p_1 \circ s_2 \circ p_2$$
and $$p_2' = 2p_2 - p_2 \circ s_2 \circ p_2
- p_2 \circ s_1 \circ p_1.$$
It follows from (1) and (2) that
$p_1' \circ s_1 = 2p_1 \circ s_1 - (2p_1 \circ s_1 -
1) - 0 = 1$.
        From (3) and (4) we find that
$p_1' \circ s_2 = 2p_1 \circ s_2 - p_1 \circ s_2 -
p_1 \circ s_2 = 0$.
Finally, using (5) we obtain that $s_1 \circ p_1' +
s_2 \circ p_2' = 1$.
Combining these results with their symmetric results
(changing the roles of $A$ and $B$), we have shown that
$(C,s_1,s_2,p_1',p_2')$ is a biproduct of $A$ and $B$
in $\AAA$, which proves our assertion.
\end{proof}
\end{proposition}

\begin{proposition}\label{ref-A.7-147}
If $\BBB$ is additive, the same holds for $\AAA$.
\begin{proof}
This immediately follows from Proposition \ref{ref-A.2-142} and
Proposition \ref{ref-A.6-146}.
\end{proof}
\end{proposition}

\begin{definition}\label{ref-A.8-148}
A pre-additive category $\AAA$ is called \emph{coherent} if $\mmod(\AAA)$ is
abelian.
\end{definition}

\begin{proposition}\label{ref-A.9-149}
If $\BBB$ is coherent, the same holds for $\AAA$.
\begin{proof}
We may assume that $\AAA,\BBB$ are small. Then the result follows from
Proposition \ref{ref-5.4-61} and Theorem \ref{ref-6.36-103}.
\end{proof}
\end{proposition}

We mention the following intrinsic characterization of coherence, which we will
not explicitly use.

\begin{definition}\cite{Krause2}
Consider a map $f:C \lra C'$ in a pre-additive category $\ccc$. A \emph{weak
cokernel} of $f$ is a map $g:C' \lra C''$ with $g \circ f = 0$ and such that
every map $h$ with $h \circ f = 0$ factorizes as $h = h' \circ g$.
\end{definition}

In a triangulated category, the cone of a morphism is a weak
cokernel.

\begin{proposition}\cite[Lemma 1]{Krause2}
An additive category $\AAA$ is coherent if and only if
$\AAA$ has weak cokernels.
\end{proposition}

\begin{remark}\label{ref-A.12-150}
If $\ccc$ is an
abelian category with enough injectives, there is an equivalence $\ccc
\cong (\mmod(\Inj(\ccc)))^{\op}$ (see Proposition \ref{ref-6.25-92}) hence
$\Inj(\ccc)$ is coherent. For a map between injectives, we can first take its
cokernel and then a mono to an injective to obtain a weak cokernel in
$\Inj(\ccc)$.
\end{remark}

\subsection{Deformations with enough injectives}

We have seen in Theorem \ref{ref-8.17-135} that if $\ccc$ is an abelian category with
enough injectives, there is an equivalence $\Def_{\ccc} \cong
\ddef_{\inj(\ccc)}$. \FORM We will now give a different approach to this fact.
We start by characterizing ``categories of injectives of abelian categories
with enough injectives''.

\begin{proposition}\label{ref-A.13-151}
For an abelian category $\ccc$ with enough injectives, $\Inj(\ccc)$ is an
additive, coherent, Karoubian category.
\begin{proof}
$\Inj(\ccc)$ is additive since products of injectives are injective
and we already noticed in Remark \ref{ref-A.4-144} that $\Inj(\ccc)$ is
Karoubian and in Remark \ref{ref-A.12-150} that $\Inj(\ccc)$ is coherent.
\end{proof}
\end{proposition}

\begin{proposition}\label{ref-A.14-152}
For an additive, coherent, Karoubian category $\AAA$, the functor
$$\AAA\op \lra \mmod(\AAA): A \longmapsto (A,-)$$
induces an equivalence of categories
$$\AAA\op \lra \PPP$$
to the full subcategory $\PPP$ of projectives of
$\mmod(\AAA)$.
\begin{proof}
Take an object $P$ of $\PPP$ and
consider a presentation $$\oplus_{i=1}^n \AAA(A_i,-) \lra
\oplus_{j=1}^m \AAA(A_j,-) \lra P \lra 0.$$ Since $P$ is
projective and $\AAA$
is additive, $P$ is a retract of the functor
$\AAA(\oplus_{j=1}^mA_j,-)$.
But since $\AAA$ is Karoubian, it follows that $P$ is itself
representable.
\end{proof}
\end{proposition}

In other words, for an additive, coherent, Karoubian category $\AAA$, there
is an equivalence $\AAA \cong \Inj((\mmod(\AAA))^{\op})$. Combining this with
Proposition \ref{ref-A.13-151}, we have thus characterized the ``categories of
injectives of abelian categories with enough injectives'' as being precisely
the additive, coherent, Karoubian categories.

Let $\BBB$ be a coherent $S$-linear category.
Since $\mmod(\BBB) \cong \vvv-\Fp(\vvv-\Mod(\BBB))$ where $\BBB \in \vvv$, we
obtain a pseudo natural transformation
\begin{equation}
\label{ref-A.1-153}
\ddef_{\BBB} \lra \Def_{(\mmod(\BBB))^{\op}},
\end{equation}
which is actually an equivalence by Theorem \ref{ref-8.16-133}, Theorem
\ref{ref-8.8-125} and Proposition \ref{ref-8.7-121}.
Taking $\BBB = \Inj(\ddd)$, this is precisely the equivalence
(\ref{ref-8.3-136}) of section \ref{ref-8.8-134}. We will now suggest another proof
of this equivalence. \FORM

\begin{proposition}\label{ref-A.15-154}
Every flat, nilpotent abelian deformation $\ddd \lra \ccc$ with enough
injectives induces a flat linear deformation $\HomS: \Inj(\ccc) \lra
\Inj(\ddd)$.
\begin{proof}
Since taking injectives preserves equivalence of categories, we may consider
$\ccc_S \lra \ccc$.
By Proposition \ref{ref-2.9-12}(\ref{ref-5-19}) and Proposition \ref{ref-4.3-43},
the map
$$\ccc(C,E)
\lra \ccc(\Hom_R(S,C),E) \lra \ccc(\Hom_R(S,C),\Hom_R(S,E))$$
induces an isomorphism $S \otimes_R \ccc(C,E) \lra
\ccc(\Hom_R(S,C),\Hom_R(S,E))$ when $E$ is injective.
It then follows by Corollary \ref{ref-6.15-81} that $\HomS: \Inj(\ccc)^S
\lra \Inj(\ccc_S)$ is an equivalence.
The flatness of $\Inj(\ccc)$ follows from \ref{ref-2.9-12}(\ref{ref-6-20}).
\end{proof}
\end{proposition}

Let $\ddd$ be an abelian $S$-linear category with enough injectives.
Proposition
\ref{ref-A.15-154} yields a pseudo natural transformation

\begin{equation}
\label{ref-A.2-155}
\Def_{\ddd} \lra \ddef_{\Inj(\ddd)}.
\end{equation}

By Theorem \ref{ref-6.16-82},
any $\ccc \in \Def_{\ddd}$ has enough injectives, and by Propositions
\ref{ref-A.13-151}, \ref{ref-A.7-147}, \ref{ref-A.5-145} and
\ref{ref-A.9-149}, any $\AAA \in \ddef_{\BBB}$ is an additive,
coherent, Karoubian category. This allows one to prove the
following

\begin{theorem}
(\ref{ref-A.1-153}) and (\ref{ref-A.2-155}) are equivalences of deformation
pseudofunctors.
\end{theorem}

\section{\strut}
\label{ref-B-156}

In this Appendix, we consider an alternative deformation pseudo functor
$\ddef^s_{\BBB}$ that can be used to study linear deformations of an
$S$-linear category $\BBB$, and we study its relation with $\ddef_{\BBB}$.

For an $S$-linear category $\BBB$ and for $R \in \Rng^0/S$, consider the
following groupoid $\ddef^s_{\BBB}$: the objects of $\ddef^s_{\BBB}$ are
strict $R$-deformations $\AAA \lra \BBB$. The morphisms of $\ddef^s_{\BBB}$
are isomorphisms of deformations.

There are obvious functors
$$\sigma_R: \ddef^s_{\BBB}(R) \lra \ddef_{\BBB}(R).$$
constituting a pseudonatural transformation $\sigma: \ddef^s_{\BBB}
\lra \ddef_{\BBB}.$ As every functor between groupoids does, $\sigma_R$
induces a function
$$\Sk(\sigma_R): \Sk(\ddef^s_{\BBB}(R)) \lra
\Sk(\ddef_{\BBB}(R)).$$
We will show that for $R$ in $\Rng^0/S$, $\Sk(\sigma_R)$ is a bijection.

\begin{proposition}\label{ref-B.1-157}
For $R$ in $\Rng^0/S$, $\sigma_R$ is full.
\begin{proof}
Consider two strict left $R$-deformations
$f_1:
\AAA_1
\lra
\BBB$ and $f_2: \AAA_2 \lra \BBB$ and an equivalence of
deformations
$\varphi: \AAA_1 \lra \AAA_2$. We will construct an
isomorphism of deformations
$\varphi':\AAA_1 \lra \AAA_2$ that is naturally isomorphic to $\varphi$.
Consider a natural isomorphism $\eta: f_1 \lra f_2 \circ \varphi$. For every
$\AAA$-map $a: A \lra A'$, there is a commutative square
$$\xymatrix{f_1(A) \ar[d]_{f_1(a)} \ar[r]^-{\eta_A} & f_2(\varphi(A))
\ar[d]^{f_2(\varphi(a))} \\ f_1(A') \ar[r]_-{\eta_{A'}} & f_2(\varphi(A')).}$$
Since $f_2$ is a strict deformation, every $A$ in $\AAA_1$ determines a unique
object $\varphi'(A)$ in $\AAA_2$ satisfying $$f_2(\varphi'(A)) = f_1(A).$$
By Proposition \ref{ref-A.1-141}, we can lift the $\BBB$-isomorphisms $$\eta_A:
f_2(\varphi'(A)) \lra f_2(\varphi(A))$$ to $\AAA_2$-isomorphisms
$$\mu_A: \varphi'(A) \lra \varphi(A)$$ with $f_2(\mu_A) = \eta_A$. We can now
define $\varphi'(a)$ for $a: A \lra A'$ in $\AAA_1$ as the unique $\AAA_2$-map
making the following square commute:
$$\xymatrix{\varphi'(A) \ar[d]_{\varphi'(a)} \ar[r]^-{\mu_A} & \varphi(A)
\ar[d]^{\varphi(a)} \\ \varphi'(A') \ar[r]_-{\mu_A'} & \varphi(A')}$$
Clearly, $\varphi': \AAA_1 \lra \AAA_2$ is a functor and $\mu: \varphi' \lra
\varphi$ is a natural isomorphism.
If we apply $f_2$ to the second square and compare the resulting square with
the first square, we conclude that $$f_2 \circ \varphi' = f_1,$$
as we wanted.\end{proof}
\end{proposition}

\begin{remark}\label{ref-B.2-158}
$\sigma_R$ need in general not be faithfull.
Indeed, consider any algebra-deformation $f: A \lra B$ in which $A$ contains a
non-central invertible element $\xi$ with $f(\xi) = 1$. Then $\varphi: A \lra
A: a
\longmapsto
\xi^{-1}a\xi$ is an algebra-isomorphism different from but naturally
isomorphic to $1_A$ with $f \circ \varphi = f$.
\end{remark}

\begin{proposition}\label{ref-B.3-159}
For $R$ in $\Rng^0/S$, $\sigma_R$ is essentially surjective.
\begin{proof}
Consider a left $R$-deformation $f: \AAA \lra \BBB$. Take an inverse
equivalence $g: \BBB \lra S \otimes_R\AAA$ of the canonical equivalence $S
\otimes_Rf:
S \otimes_R\AAA
\lra \BBB$.
We construct a category $\CCC$ in the following way: the objects of $\CCC$ are
precisely the objects of $\BBB$. For objects $B$, $B'$ of $\BBB$, $\CCC(B,B')$
is defined to equal $\AAA(g(B),g(B'))$ and the composition in $\CCC$ is the
composition in $\AAA$. There is an obvious functor $\varphi: \CCC \lra \AAA$
mapping $B$ to $g(B)$ and $a:g(B) \lra g(B')$ to $a:g(B) \lra g(B')$.
$\varphi$ is clearly fully faithful. For $A$ in $\AAA$, $\varphi(f(A)) =
g(f(A))$ is isomorphic to $A$ in $S \otimes_R\AAA$. But by Proposition
\ref{ref-A.1-141}, they remain isomorphic in $\AAA$. It follows that
$\varphi$ is essentially surjective and hence an equivalence of categories.
Take a natural isomorphism $\eta: 1_{\BBB} \lra S \otimes_Rf \circ g$.
We define a functor $h: \CCC \lra \BBB$ by putting $h(B) = B$ and mapping
$a: g(B) \lra g(B')$ to the unique map $h(a)$ making the following diagram
commute:
$$\xymatrix{B \ar[r]^-{\eta_B} \ar[d]_{h(a)} & f(g(B)) \ar[d]^{f(a)} \\
B' \ar[r]_-{\eta_{B'}} & f(g(B'))}$$
It follows that the $\eta$ define a natural isomorphism $h \lra f \circ
\varphi$.
Finally, we have to show that $S \otimes_Rh: S \otimes_R\CCC \lra \BBB$ is an
isomorphism. Since
$S \otimes_Rh$ is clearly bijective on objects, it suffices that $S
\otimes_Rh$ or equivalently
$S \otimes_R(f \circ \varphi)$ is an equivalence, which is obvious.
\end{proof}
\end{proposition}

\begin{theorem}\label{ref-B.4-160}
For $R$ in $\Rng^0/S$, $\Sk(\sigma_R)$ is a bijection.
\begin{proof}
$\Sk(\sigma_R)$ is injective by \ref{ref-B.1-157} and surjective by
\ref{ref-B.3-159}.
\end{proof}
\end{theorem}

%\bibliography{defpap}

\begin{thebibliography}{10}

\bibitem{AHS}
J.~Ad{\'a}mek, H.~Herrlich, and G.~E. Strecker,{\em Abstract and concrete
categories}, Pure and Applied Mathematics,
   John Wiley \& Sons Inc., New York, 1990, The joy of cats, A
   Wiley-Interscience Publication.

\bibitem{SGA4}
M.~Artin, A.~Grothendieck, and J.~L. Verdier, {\em Theorie des topos et
   cohomologie \'etale des sch\'emas, {SGA4}, {T}ome 3}, Lecture Notes in
   Mathematics, vol. 305, Springer Verlag, 1973.

\bibitem{AZ2}
M.~Artin and J.~J. Zhang, {\em Abstract {H}ilbert schemes}, Algebr. Represent.
   Theory {\bf 4} (2001), no.~4, 305--394.

\bibitem{HCA2}
F.~Borceux, {\em Handbook of categorical algebra. 2}, Encyclopedia of
   Mathematics and its Applications, vol.~51, Cambridge University Press,
   Cambridge, 1994, Categories and structures.

\bibitem{Bourbakistructure}
N.~Bourbaki,{\em{E}l\'ements de math\'ematique. 22. {P}remi\`ere partie:
            {L}es structures fondamentales de l'analyse. {L}ivre 1:
           {T}h\'eorie des ensembles. {C}hapitre 4: {S}tructures},
Actualit\'es Sci. Ind. no. 1258, Hermann, Paris, 1957.

\bibitem{GS1}
M.~Gerstenhaber and S.~D. Schack, {\em On the deformation of algebra morphisms
   and diagrams}, Trans. Amer. Math. Soc. {\bf 279} (1983), no.~1, 1--50.

\bibitem{GS}
\bysame, {\em Algebraic cohomology and deformation theory}, Deformation theory
   of algebras and structures and applications (Il Ciocco, 1986) (Dordrecht),
   NATO Adv. Sci. Inst. Ser. C Math. Phys. Sci., vol. 247, Kluwer Acad. Publ.,
   Dordrecht, 1988, pp.~11--264.

\bibitem{GS2}
\bysame, {\em The cohomology of presheaves of algebras. {I}. {P}resheaves over
   a partially ordered set}, Trans. Amer. Math. Soc. {\bf 310} (1988), no.~1,
   135--165.

\bibitem{GerstI}
M.~Gerstenhaber, {\em On the deformation of rings and algebras}, Ann. of Math.
   (2) {\bf 79} (1964), 59--103.

\bibitem{GerstII}
\bysame, {\em On the deformation of rings and algebras. {II}}, Ann. of Math.
   {\bf 84} (1966), 1--19.

\bibitem{Groth1}
A.~Grothendieck, {\em Sur quelques points d'alg\`ebre homologiques}, T{\^o}hoku
   Math. J. (2) {\bf 9} (1957), 119--221.

\bibitem{Illusie}
L.~Illusie, {\em Existence de r\'esolutions globales}, SGA6, Lecture Notes in
   Math., vol. 225, Springer Verlag, 1971.

\bibitem{Krause1}
H.~Krause, {\em The spectrum of a module category}, Mem. Amer. Math. Soc. {\bf
   149} (2001), no.~707, x+125.

\bibitem{Krause2}
\bysame, {\em A {B}rown representability theorem via coherent functors},
   Topology {\bf 41} (2002), no.~4, 853--861.

\bibitem{lowen1}
W.~T. Lowen, {\em A generalization of the {G}abriel-{P}opescu theorem}, to
   appear.

\bibitem{low2}
\bysame, {\em An obstruction theory for objects in abelian categories}, in
   preparation.

\bibitem{LVdBII}
W.~T. Lowen and M.~Van~den Bergh, {\em Deformation theory for abelian
   categories {II}}, in preparation.

\bibitem{maclane}
S.~MacLane, {\em Categories for the working mathematician}, Springer Verlag,
   Berlin, 1971.

\bibitem{Mi}
B.~Mitchell, {\em Rings with several objects}, Advances in Math. {\bf 8}
   (1972), 1--161.

\bibitem{Neemanboek}
A.~Neeman, {\em Triangulated categories}, Annals of Mathematics Studies, vol.
   148, Princeton University Press, Princeton, NJ, 2001.

\bibitem{GP}
N.~Popesco and P.~Gabriel, {\em Caract\'erisation des cat\'egories ab\'eliennes
   avec g\'en\'erateurs et limites inductives exactes}, C. R. Acad. Sci. Paris
   {\bf 258} (1964), 4188--4190.

\bibitem{Popescu}
N.~Popescu, {\em Abelian categories with applications to rings and modules},
   Academic Press, London, 1973, London Mathematical Society Monographs, No. 3.

\bibitem{Ta}
L.~A. Takhtadjian, {\em Noncommutive homology of quantum tori}, Functional
   Anal. Appl. {\bf 23} (1989), 147--149.

\end{thebibliography}
%\bibliographystyle{amsabbrv}
\def\cprime{$'$}
\ifx\undefined\bysame
\newcommand{\bysame}{\leavevmode\hbox to3em{\hrulefill}\,}
\fi

\end{document}